\documentclass[onefignum,onetabnum]{siamonline171218}
\usepackage{color}
\usepackage{hyperref}       % hyperlinks
\usepackage{url}            % simple URL typesetting
\usepackage{makecell}
\usepackage{amsmath,amsfonts,amssymb}
\usepackage{algpseudocode}
\usepackage{booktabs}       % professional-quality tables
\usepackage{nicefrac}       % compact symbols for 1/2, etc.
\usepackage{microtype}      % microtypography
\usepackage{lipsum}
\usepackage{graphicx}
\usepackage{ragged2e}
\usepackage{multirow}
\usepackage{lineno}
\usepackage{caption}
\usepackage{pifont}
\headers{Blood Flow in Elastic Vessel Using PINN}{H. Zhang, R. Chan and X.C. Tai}

\title{A Meshless Solver for Blood Flow Simulations in Elastic Vessels Using Physics-Informed Neural Network
\thanks{Submitted to the editors DATE.
}}

\author{Han Zhang \thanks{Department of Mathematics, City University of Hong Kong and Hong Kong Centre for Cerebro-Cardiovascular Health EngineeringCity University of Hong Kong, Hong Kong, China. (\email{hzhang863-c@my.cityu.edu.hk})}
\and Raymond Chan \thanks{Department of Mathematics, City University of Hong Kong and Hong Kong Centre for Cerebro-Cardiovascular Health EngineeringCity University of Hong Kong, Hong Kong, China. (\email{raymond.chan@cityu.edu.hk})}
\and Xue-Cheng Tai \thanks{Norwegian Research Centre (NORCE), Nygardsgaten 112, 5008 Bergen, Norway. (\email{xtai@norceresearch.no}, \email{xuechengtai@gmail.com})}
}

\begin{document}
\maketitle
\begin{abstract}

Investigating blood flow in the cardiovascular system is crucial for assessing cardiovascular health. Computational approaches offer some non-invasive alternatives to measure blood flow dynamics. Numerical simulations based on traditional methods such as finite-element and other numerical discretizations have been extensively studied and have yielded excellent results. However, adapting these methods to real-life simulations remains a complex task.
In this paper, we propose a method that offers flexibility and can efficiently handle real-life simulations. We suggest utilizing the physics-informed neural network (PINN) to solve the Navier-Stokes equation in a deformable domain, specifically addressing the simulation of blood flow in elastic vessels. Our approach models blood flow using an incompressible, viscous Navier-Stokes equation in an Arbitrary Lagrangian-Eulerian form. The mechanical model for the vessel wall structure is formulated by an equation of Newton's second law of momentum and linear elasticity to the force exerted by the fluid flow.
Our method is a mesh-free approach that eliminates the need for discretization and meshing of the computational domain. This makes it highly efficient in solving simulations involving complex geometries. Additionally, with the availability of well-developed open-source machine learning framework packages and parallel modules, our method can easily be accelerated through GPU computing and parallel computing.
To evaluate our approach, we conducted experiments on regular cylinder vessels as well as vessels with plaque on their walls. We compared our results to a solution calculated by Finite Element Methods using a dense grid and small time steps, which we considered as the ground truth solution. We report the relative error and the time consumed to solve the problem, highlighting the advantages of our method.

\end{abstract}

\begin{keywords}
  Fluid-structure interaction, Physics-informed neural network, blood flow simulation, Arbitrary Lagrangian-Eulerian, Computational Fluid Dynamics
\end{keywords}
\section{Introduction}

The measurement of blood velocity and pressure is a fundamental aspect of medical practice, providing crucial insights into the dynamics of the cardiovascular system. Despite the widespread application of established techniques such as Fractional Flow Reserve (FFR) \cite{pijls1996measurement}, it is important to recognize their limitations, particularly in specific clinical scenarios \cite{toth2014revascularization,dattilo2012contemporary}. For instance, the invasive operation of FFR assessment can increase patient risks, potentially leading to complications or exacerbating existing medical conditions. Inserting an intravascular wire into the vessel can disrupt the original flow dynamics by altering the spatial configuration of the inner lumen \cite{yan2022impact}. To address these challenges, research has been directed towards using computational approaches, which involves geometry reconstruction from medical images \cite{papafaklis2014fast,norgaard2014diagnostic} and numerically solving for a properly modeled Fluid-Structure Interaction (FSI) problem \cite{moore1999accuracy}.

The FSI problem is modeled as a coupled problem of the blood flow and the vessel structures. Specifically, the fluid stress produces the deformation of structures and will reversely change the fluid dynamics. To describe and solve the blood flow simulation problem, the Arbitrary Lagrangian-Eulerian (ALE) formulation is a popular choice. Within this formulation, the displacement of the computational domain is introduced as a variable in the coupled system, ensuring the satisfaction of coupling conditions on the fluid-structure interface. While previous research has used linear elasticity to model the vessel wall \cite{quarteroni2004mathematical}, others have modeled it as poroelastic media by coupling the Navier-Stokes equation with the Biot system \cite{badia2009coupling}. Quaini et al. \cite{quaini2012validation} validated the simulated results, which use a semi-implicit, monolithic method, through a mock heart chamber using an elastic aperture.

A variety of methods have been developed to solve these models. Fernandez et al. \cite{fernandez2007projection} introduced a semi-implicit coupling method that implicitly couples the pressure stress to the structure while treating the nonlinearity due to convection and geometrical nonlinearities explicitly. To further reduce computational costs, Badia et al. \cite{badia2008splitting} designed two methods based on the inexact factorization of the linearized fluid-structure system. Barker and Cai \cite{barker2010scalable} developed a method that couples the fluid to the structure monolithically, making it robust to changes in physical parameters and large deformation of the fluid domain. A loosely coupled scheme that is unconditionally stable was proposed by Bukavc et al. \cite{bukavc2013fluid}, based on a newly modified Lie operator splitting that decouples the fluid and structure sub-problems. Other approaches include the space-time formulation \cite{bazilevs2008isogeometric,tezduyar2007modelling}, the higher-order discontinuous Galerkin method \cite{wang2018higher}, the immersed boundary method \cite{wang2008immersed,ge2007numerical}, and the coupled momentum method \cite{figueroa2006coupled,zhou2010cardiovascular}.

However, numerically solving such an FSI problem remains to be challenging. The mesh must be carefully discretized since the deformation of the mesh over time will change the mesh geometry, resulting in the failure of convergence during the solving process \cite{wick2011fluid,sun2022coupled}. Current autonomous software like Gmsh \cite{geuzaine2009gmsh} can produce sub-optimal meshes that can cause methods like the Finite Element Method (FEM) to fail in complex FSI problems. To address the issues caused by large structural displacement, Basting et al. \cite{basting2017extended} proposed the Extended ALE Method, which relies on a variational mesh optimization technique where mesh alignment with the structure is achieved via a constraint. However, in many cases, manual corrections on the mesh are required to make a result. These corrections may be needed multiple times, with validation required after each correction. Such exhausting operations make solving problems on complex geometry hard and expensive. Besides, solutions by conventional methods are not continuously defined in the domain of our model, leading to further errors if an interpolation is needed to find a solution at a particular point that is not a node in the discretized grid. What's worse, to acquire an accurate solution, the conventional mesh-based methods require a high-resolution configuration of dense grids and small time steps. Such settings would dramatically increase the computation workload and memory consumption, slowing and preventing a fast and satisfying result acquisition.

To address these challenges, we propose a neural solver for the Navier-Stokes equation in deformable vessels using Physics-Informed Neural Networks (PINNs)\cite{raissi2019physics}. In our methods, we take some networks to produce displacement, pressure and velocity of our problem. To acquire those networks, we first convert the partial differential equations of the formulated physical model into residual forms and use them as loss functions. Next, the automatic differentiation is used to represent all the differential operators. Then, we train the networks to minimize those losses. With those losses getting approximately zero, the networks can be an approximation of the functions for displacement, pressure, and velocity that are used to solve the modeled problem. 

Benefit from the PINNs method, the proposed method is a mesh-free approach, eliminating the need for discretization and meshing for the computational domain. This feature makes it highly efficient and straightforward to add other objects into the geometry (e.g., blood plaque) for simulation. To evaluate the performance of our methods, we performed several experiments. The first experiment simulates the blood flow in an elastic cylinder-like vessel. For comparison, we use Finite Element Methods to solve the same problem using different grid resolutions and time steps. The results revealed that our approach can produce comparable results to FEM in high-resolution grids while consuming notably less time. Furthermore, we perform the blood flow simulation on a vessel with a plaque attained further to illustrate the robustness of our methods on complex geometry. In addition, we performed self-ablation on networks with different depths and widths to show how we should choose the architecture for the neural networks for our task. Since our methods are easy to parallelize, an evaluation of how our method can be accelerated by using multiple GPUs is performed. Last but not least, the mechanical analysis of how different sizes of plaque might affect the blood flow dynamics is done by using our methods.

Our contributions are summarized as follows: 
\begin{enumerate} 
    \item We incorporated the Navier-Stokes equation in the Arbitrary Lagrangian-Eulerian (ALE) form into the Physics-Informed Neural Network to tackle fluid-structure interaction problems in a meshless approach.
    \item We designed a sequential-alternative training approach by conceptualizing our problem as a multi-objective optimization task encompassing fluid dynamics and deformation. This approach allows for stable training and better convergence by focusing on optimizing subsets of losses sequentially.
    \item We introduced a novel architecture featuring separate networks for predicting velocity and pressure. This design choice aligns with the distinct operators governing these variables in the Navier-Stokes equation, thereby enabling effective capture of their unique characteristics and complexities.
    \item We proposed a new architecture utilizing alternative activation functions by \textit{ReLU} and \textit{Sigmoid}. This modification results in faster training and improved convergence dynamics. 
\end{enumerate}

% \begin{algorithm}
% \caption{Training Scheme for Fluid Problem}
% \label{alg:optimizefluid1}
% \begin{algorithmic}[1]
% \State $k=1$;
% \For{$k=1:\frac{m_{\text{fluid}}}{(m1+m2)}$}
%     \State $\boldsymbol{\theta}_{u} \gets $ Optimize fluid problem using Eq. 1 for $m_1$ epochs;
%     \State $\boldsymbol{\theta}_{p} \gets $ Optimize fluid problem using Eq. 2 for $m_2$ epochs;
% \EndFor
% \end{algorithmic}
% \end{algorithm}
% \begin{algorithm}
% \caption{Fluid-Structure Interaction Training Scheme}
% \label{alg:trainingscheme1}
% \begin{algorithmic}[1]
% \State Initialize $N_d$ with \textit{zero\_initialization};
% \State $\alpha_{\text{ns}} \gets 0$;
% \State $\boldsymbol{\theta}_{u},\boldsymbol{\theta}_{p} \gets $ Optimize fluid problem using Alg. \ref{alg:optimizefluid1} for $m_{\text{fluid}}$ epochs;
% \State $\alpha_{\text{ns}} \gets 10^{-8}$;
% \For{i=1:5}
%     \State $\alpha_{\text{ns}} \gets 10\times\alpha_{\text{ns}}$;
%     \State $\boldsymbol{\theta}_{u},\boldsymbol{\theta}_{p} \gets $ Optimize fluid problem using Alg. \ref{alg:optimizefluid1} for $m_{\text{fluid}}$ epochs;
% \EndFor
% \For{i=1:F \textbf{and} $\mathcal{L}_{\text{fluid}}$ doesn't converge}
%     \State $\boldsymbol{\theta}_{d} \gets $ Optimize solid problem using Eq. 3 for $m_{\text{solid}}$ epochs;
%     \State $\boldsymbol{\theta}_{u},\boldsymbol{\theta}_{p} \gets $ Optimize fluid problem using Alg. \ref{alg:optimizefluid1} for $m_{\text{fluid}}$ epochs;
% \EndFor
% \end{algorithmic}
% \end{algorithm}
\section{Problem Formulation}\label{sec:problem}
\begin{figure}[!t]
    \centering
    \includegraphics[width=\textwidth]{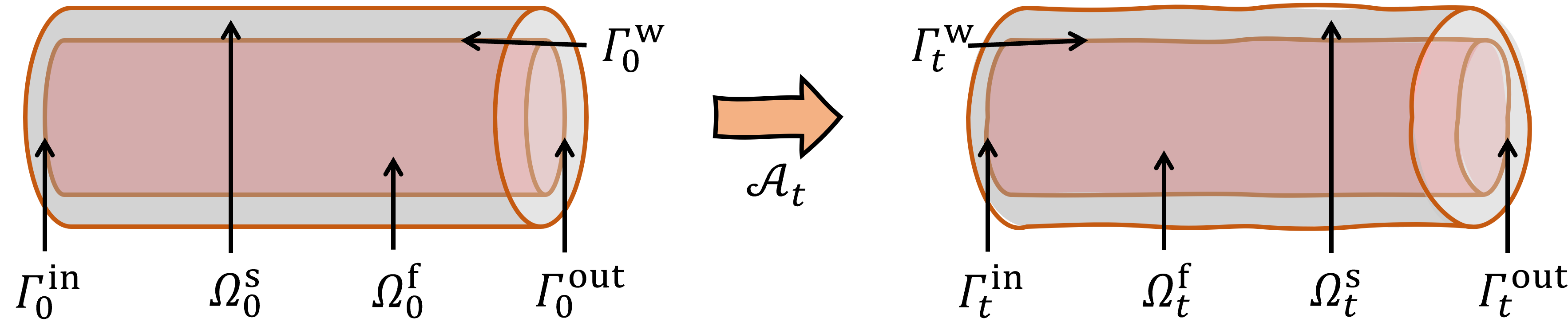}
    \caption{The illustration for the Arbitrary Lagrangian Euler form. Left: reference configuration; Right: current configuration.}
    \label{fig:ALEmapping}
\end{figure}

In this section, we introduce the mathematical model for our problem, which is to simulate the flow of blood in a deformable vessel. In essence, the blood flow, assumed to be incompressible, is described by the Navier-Stokes equations in the pressure-velocity formulation. Since the computational domain is moving, we employ the Arbitrary Lagrangian Euler (ALE) form of the Navier-Stokes equations in this study. The vessel wall, assumed to be linear-elastic, is continuously moving due to the force exerted by the fluid. With the movement of the structures, the computational domain for the fluid is also moving and is computed by the harmonic extension from the displacement of the structure.

To accurately represent the deformation and simulate the flow, it is crucial to properly model the movement of the computational domain, as well as the fluid velocity and pressure. As illustrated in Figure \ref{fig:ALEmapping}, we start with the reference configuration $\Omega_0$, which serves as the parameterization for the entire domain required for our computations. Dividing this region into two parts, we have the fluid domain $\Omega_0^{\text{f}}$, through which the flow occurs, and the structural domain $\Omega_0^{\text{s}}$, representing the vessel wall. The interface between the fluid and structural domains is denoted as $\Gamma_0^{\text{w}}$. It's important to note that most of the interaction takes place at this interface. The inlet boundary is labeled as $\Gamma_0^{\text{in}}$ and the outlet boundary is as $\Gamma_0^{\text{out}}$. Each point in the domain of the reference configuration is denoted as the problem coordinates $\boldsymbol{x}_0 = (x_0, y_0, z_0)^\intercal$, where the bold symbols are used throughout this paper to signify matrices or column vectors. As the domain continues to move, we use $\Omega_t$ to represent the current configuration at time $t$. The symbols $\Omega_t^{\text{f}}$, $\Omega_t^{\text{s}}$, $\Gamma_0^{\text{w}}$, $\Gamma_0^{\text{in}}$, and $\Gamma_0^{\text{out}}$ denote the fluid domain, the structural domain, the interface, the inlet, and the outlet, respectively. The coordinates of points in the current configuration are referred to as the physical coordinates $(x_t, y_t, z_t)^\intercal$.

Following Quarteroni et al. \cite[P. 73, Sec. 18]{quarteroni2004mathematical}, we introduce the ALE mapping defined as:
\begin{equation}
\mathcal{A}_{t}: \Omega_{0} \rightarrow \Omega_{t}, \quad \boldsymbol{x}_0 \rightarrow \boldsymbol{x}_t = \boldsymbol{x}_0 + \boldsymbol{d}(\boldsymbol{x}_0, t) = \mathcal{A}_{t}(\boldsymbol{x}_0)\quad \text{ in } \Omega_t^{\text{f}},
\label{eq:ALEmapping}
\end{equation}
where the displacement function $\boldsymbol{d}(\cdot)$ characterizes the domain movement within the reference configuration, and $\Tilde{\boldsymbol{w}} = \frac{\partial \boldsymbol{d}}{\partial t}$ denotes the domain velocity within $\Omega_0$. In the current domain $\Omega_t$, the domain velocity is defined as $\boldsymbol{w} = \Tilde{\boldsymbol{w}} \circ \mathcal{A}_t^{-1}$.

Suppose no external force is applied to the flow, the incompressible unsteady Navier-Stokes equation which describes the movement of a viscous fluid in a moving domain \cite[P. 74, Sec. 18]{quarteroni2004mathematical} is:
\begin{equation}
\begin{aligned}
\rho_\text{f} \frac{D}{D t} \boldsymbol{u}+ \rho_\text{f} [(\boldsymbol{u}-\boldsymbol{w}) \cdot \nabla] \boldsymbol{u} -  \nabla \cdot \boldsymbol{\sigma}_\text{f}=\boldsymbol{0} \\
\nabla\cdot \boldsymbol{u}=0 
\label{eq:ns}
\end{aligned}\quad \text{ in } \Omega_t^{\text{f}}, \ t \in [0,T],
\end{equation}
where $\rho_\text{f}$ is the fluid density, $\boldsymbol{u}$ is the fluid velocity, $\boldsymbol{\sigma}_\text{f} = -P\boldsymbol{I} + 2\mu \frac{ \nabla \boldsymbol{u} + (\nabla \boldsymbol{u})^\intercal}{2}$ is the Cauchy stress tensor, $P$ is fluid pressure, $\mu$ is the viscosity of the fluid and is a positive quantity. The term $\frac{D}{D t} \boldsymbol{u} = \frac{\partial \boldsymbol{u}}{\partial t} + \boldsymbol{w} \cdot \nabla \boldsymbol{u}$ represents the ALE derivative and $\nabla$ is the gradient with respect to $\boldsymbol{x}_t$ within $\Omega_t^{\text{f}}$. Here, $(\nabla \boldsymbol{u})^\intercal$ is the transpose of $\nabla \boldsymbol{u}$.

In modeling the structure of the vessel, we make several simplifying assumptions: (1) The thickness of the vessel wall is small enough to permit a shell-type representation of the vessel geometry; (2) The reference vessel configuration is approximated as a circular cylindrical surface with straight axes, and the displacements occur primarily in the radial direction; (3) The mechanical properties of the vessel wall follow linear elasticity. While these assumptions simplify our problem, it's important to note that extending the model to more complex scenarios is a straightforward possibility.

% \begin{figure}[!t]
%     \centering
%     \includegraphics[width=0.8\textwidth]{img/CoordinateSystem.png}
%     \caption{The coordinate system illustration. A 3D geometry can be represented using a Cartesian coordinate system as $(x,y,z)$ (left). A cylindrical model of the vessel geometry may be described by using a curvilinear cylindrical coordinate system $r = R(\theta, z;t)$ with the corresponding base unit vectors $\mathbf{e}_r$, $\mathbf{e}_\theta$ and $\mathbf{e}_z$, where $\mathbf{e}_z$ is aligned with the axis of the artery (right).}
%     \label{fig:coordinatesystem}
% \end{figure}

Assuming the vessel has a cylinder-like geometry at time $0$, which is also the reference domain, the vessel wall \cite[P. 60, Sec. 16]{quarteroni2004mathematical} can be rewritten as:
\begin{equation}
    \Gamma_{0}^{\text{w}}=\left\{(r, \theta, z): r=R_{0}(\theta,z), \theta \in[0,2 \pi), z \in[0, L]\right\},
\end{equation}
where $R_{0}(\theta,z)$ is the radius for the vessel lumen at $(\theta,z)$. Assume the vessel geometry in reference configuration to be a cylinder-like geometry, we have $R_{0}(\theta,z)=r_0$ with $r_0 \in \mathbb{R}$ is a constant, representing the initial radius of the cylindrical vessel geometry, and $L$ denotes the length of the arterial element under consideration. In our cylindrical coordinate system $(r, \theta, z)$, we assume that displacements are mainly in the radial direction and other components are negligible. Thus, we can simplify the displacement $\boldsymbol{d}$ to have only the radial component ${\eta}$, such that:
\begin{equation}
    \boldsymbol{d}(r, \theta, z; t) = \eta(\theta, z;t) \boldsymbol{e}_{r},
    \label{eq:eta}
\end{equation}
Here, $\eta(\theta, z;t)$ is a function for the displacement at the vessel wall, and $\boldsymbol{e}_{r}$ is the unit vector along the radial direction.

Given these settings, we can write the domain for the vessel wall $\Gamma_{t}^{w}$ at time $t$ as:
\begin{equation}    
    \Gamma_{t}^{\text{w}}=\{(r, \theta, z): r=R(\theta, z ; t), \theta \in[0,2 \pi), z \in[0, L]\},
\end{equation}
where $R(\theta, z ; t) = \eta(\theta,z; t) \boldsymbol{e}_r + R_0(\theta, z)$.

The same as \cite[P.64, Sec. 16]{quarteroni2004mathematical}, based on linear elasticity and Newton's second law of momentum, we equate the internal force of the structure to the force exerted from the fluid, obtaining the below stress continuity equation on the vessel wall:
\begin{equation}
% \rho_{w} h_{0} R_{0} \frac{\partial^{2} \eta}{\partial t^{2}}\mathrm{d} \theta \mathrm{d} z + \frac{E h_{0}}{1-\xi^{2}} \frac{\eta}{R}\mathrm{d} \theta \mathrm{d} z=-(2 \mu \boldsymbol{D}(\boldsymbol{u}) \cdot \boldsymbol{n}) \cdot \boldsymbol{e}_{r} g R_{0}\mathrm{d} \theta \mathrm{d} z + p R \mathrm{d} \theta \mathrm{d} z \quad\text{ on }\Gamma_{0}^{\text{w}},
\rho_{w} h_{0} R_{0} \frac{\partial^{2} \eta}{\partial t^{2}} + \frac{E h_{0}}{1-\xi^{2}} \frac{\eta}{R}=-(2 \mu \boldsymbol{D}(\boldsymbol{u}) \cdot \boldsymbol{n}) \cdot \boldsymbol{e}_{r} g R_{0} + P R  \quad\text{ on }\Gamma_{0}^{\text{w}},
\label{eq:CompleteStressContinuity}
\end{equation}
where $g=\frac{R}{R_{0}} \sqrt{1+\left(\frac{1}{R} \frac{\partial R}{\partial \theta}\right)^{2}+\left(\frac{\partial R}{\partial z}\right)^{2}}$, $\boldsymbol{D}(\boldsymbol{u}) = \frac{ \nabla \boldsymbol{u} + (\nabla \boldsymbol{u})^\intercal}{2}$, $\rho_\text{w}$ is the vessel wall density, $h_0$ is the vessel wall thickness, $\xi$ is the Poisson ratio and $E$ is the Young modulus\cite{quarteroni2004mathematical}.

Regard $\boldsymbol{u}$ as a known, divide \eqref{eq:CompleteStressContinuity} by $\rho_{w} h_{0} R_{0}$, we obtain the independent ring model for the stress continuity equation \cite[Eq. (16.12), P. 60, Sec. 16]{quarteroni2004mathematical}:
\begin{equation}
\frac{\partial^2 \eta}{\partial t^2}+b \eta = H \quad \text{ in } \Gamma_0^{\text{w}},
\label{eq:StressContinuity}
\end{equation}
where $b = \frac{E}{\rho_\text{s} (1 - \xi^2 ) R_0^2}$, and $H = \frac{1}{\rho_\text{s} h_0} \left[ \frac{R}{R_0} P - g \mu \left((\nabla\boldsymbol{u} + (\nabla \boldsymbol{u})^\intercal ) \cdot \boldsymbol{n}\right) \cdot \boldsymbol{e}_r \right]$.

The structural domain in the reference configuration is then denoted as:
\begin{equation}
    \Omega_{0}^{s}=\left\{(r, \theta, z): r \in [r_0,r_0+h_0], \theta \in[0,2 \pi), z \in[0, L]\right\}.
\end{equation}

We assume that the magnitude of the displacement remains constant along the radial direction. Then, the displacement $\boldsymbol{d}$ in the structure domain on the reference configuration is expressed as
\begin{equation}
\boldsymbol{d}(r,\theta,z;t) = \eta(\theta,z;t) \boldsymbol{e}_r + r  \boldsymbol{e}_r \quad\text{ for } (r,\theta,z)\in \Omega_{0}^{s}, t\in [0,T].
\end{equation}

In addition to modeling, we need to define the initial and boundary conditions. Given the variety of specific problem settings, we provide a general representation and define them specifically in the experimental section. 

For Navier Stokes equation, which we can refer to as the fluid problem, the boundary condition is given by:
\begin{equation}
    F^{\text{f}}_{\text{bdr}}(\boldsymbol{u},P) = 0 \quad\text{ on } \Gamma_t, t\in [0,T],
    \label{eq:fluidboundary}
\end{equation}
where $\Gamma_t = \Gamma^{\text{in}}_t \cup \Gamma^{\text{out}}_t \cup \Gamma^{\text{w}}_t$. Note that \eqref{eq:fluidboundary} includes the inlet boundary condition, outlet boundary condition and the interface condition, which will be explicitly given later in Section.\ref{sec:experiment}.

The initial condition can be written as :
\begin{equation}
    F^{\text{f}}_{\text{init}}(\boldsymbol{u},P) = 0 \quad\text{ in } \Omega_0^{\text{f}}, t = 0.
    \label{eq:fluidinitial}
\end{equation}

Similarly, we write the boundary condition and initial condition for the stress continuity equation \eqref{eq:StressContinuity} as:
\begin{equation}
    F^{\text{s}}_{\text{bdr}}(\boldsymbol{d}) = 0 \quad\text{ on } \partial \Gamma^{\text{w}}_0, t\in [0,T],
    \label{eq:structureboundary}
\end{equation}
and,
\begin{equation}
    F^{\text{s}}_{\text{init}}(\boldsymbol{d}) = 0 \quad\text{ on } \Gamma^{\text{w}}_0,t=0.
    \label{eq:structureinitial}
\end{equation}

With the deformation of the structure domain modeled, we model the deformation of the fluid domain as the harmonic extension from the deformation of the structure, which well-defined the computation domain:
\begin{equation}
    \Delta \boldsymbol{d} = \boldsymbol{0} \quad \text{ in }\Omega_0^{\text{f}},t\in[0,T].
    \label{eq:harmonicextension}
\end{equation}

We may then recognize the coupling between the fluid and the structure models are interactive. The solution to the fluid problem provides stress to the structure, thus the value of $H$ in \eqref{eq:StressContinuity}. The movement of the vessel wall changes the geometry on which the fluid equations must be solved. In addition, the boundary conditions for the fluid velocity in correspondence to the vessel wall are no longer homogeneous Dirichlet conditions, but they impose equality between the fluid and the structure velocity. They express the fact that the fluid particle in correspondence with the vessel wall should move at the same velocity as the wall \cite{quarteroni2004mathematical}.
\section{Our Proposed Numerical Methods}\label{sec:method}
In this section, we will introduce the fundamental concepts of neural networks (NNs) and elucidate how NNs can be employed to approximate the solution of our proposed physical model.

\subsection{Fully connected neural network}

Denoting the set of positive integers as $\mathbb{N}^{+}$, we define a simple nonlinear function, labeled as $h_{\ell}: \mathbb{R}^{M_{\ell-1}} \rightarrow \mathbb{R}^{M_{\ell}}$, where $M_{\ell} \in \mathbb{N}^{+}$ for $\ell=0, 1, \ldots, L$. The expression for $h_{\ell}$ can be represented as 
\begin{equation}
    h_{\ell}\left(\boldsymbol{x}_{\ell}\right):=\phi\left(\boldsymbol{W}_{\ell} \boldsymbol{x}_{\ell}+\boldsymbol{b}_{\ell}\right),    
\end{equation}
where $\boldsymbol{W}_{\ell} \in \mathbb{R}^{M_{\ell} \times M_{\ell-1}}$ is the weight parameter; $\boldsymbol{b}_{\ell} \in \mathbb{R}^{M_{\ell}}$ is the bias parameter; the function $\phi(\boldsymbol{y})$ represents a specified acitivation function applied element-wise to a vector $\boldsymbol{y}$ and output another vector of the same size. Commonly used activation functions include the rectified linear unit $\phi_{\text{r}}=\max{(0, y)}$ and the sigmoidal function $\phi_{\text{s}}=\left(1+e^{-y}\right)^{-1}$.

Set the input dimension $M_{0}=d_{\text{in}}$, the output dimension $M_{L} = d_{\text{out}}$. Then an FNN $N: \mathbb{R}^{d_{\text{in}}} \rightarrow \mathbb{R}^{d_{\text{out}}}$ is formulated as the composition of these $L$ simple nonlinear functions. Then, for an input $\boldsymbol{x}_0=(x_1,x_2,\cdots,x_{d_{\text{in}}-1})^\intercal$ and temporal coordinate $t$ into:
\begin{equation}
    \hat{\boldsymbol{x}} = (x_1,x_2,\cdots,x_{d-1},t)^\intercal,
\end{equation}
we have
\begin{equation}
    N(\boldsymbol{x},t ; {\boldsymbol{\theta}})=h_{L} \circ h_{L-1} \circ h_{L-2} \circ \cdots \circ h_{1}(\hat{\boldsymbol{x}}) \quad \text { for } \hat{\boldsymbol{x}} \in \mathbb{R}^{d},
\end{equation}
where ${\boldsymbol{\theta}}:=\left\{\boldsymbol{W}_{\ell}, \boldsymbol{b}_{\ell}: 1 \leq \ell \leq L\right\}$, denotes the set of all weighting parameters that defined the neural network. Here $M_{\ell}$ is called the width of the $\ell$-th layer, and $L$ is called the depth. These parameters collectively define the architecture of a fully connected neural network (FNN). By specifying the activation function $\phi$, the depth $L$, and the sequence of widths $\{ M_{\ell} \}_{\ell=1}^{L}$, the architecture of the FNN is determined. However, the parameters within ${\boldsymbol{\theta}}$ remain undetermined and need to be initialized and optimized. For simplicity, we will focus on an architecture where the width is consistent across all layers, specifically $M_{\ell}=M$ for $\ell$ ranging from $1$ to $L-1$. We use $\mathcal{F}_{L, M, \sigma}$ to denote the set of all FNNs characterized by a depth of $L$, a uniform width of $M$, and an activation function $\sigma$ \cite{gu2023deep}. In this sense, any two functions within $\mathcal{F}_{L, M, \sigma}$ only differ in the values of their respective weighting parameters ${\boldsymbol{\theta}}$.

% The total count of weighting parameters within the network can be determined. The input layer with $\ell=1$ contains $d_{\text{in}} M + M$ parameters, the output layer has $d_{\text{out}}M + d_{\text{out}}$ parameters, and the hidden layers with $\ell=2, \ldots, L-1$ comprise a combined total of $(L-2) (M^2 + M)$ parameters. Therefore, the cardinality of $\theta$ is calculated as $|{\boldsymbol{\theta}}|=d_{\text{in}} M + M + d_{\text{out}}M + d_{\text{out}} + (L-2) (M^2 + M) = (L-2) M^{2} + (d_{\text{in}}+d_{\text{out}}+L-1) M + d_{\text{out}}$.

% (d -2)*(w1**2) + (3+2+ d -1)*w1+2 = 3200 + 282 = 3482
% (d -2)*(w2**2) + (3+1+ d -1)*w2+1 = 800 + 131 = 931

To ensure the smoothness of the solutions approximated by neural networks, the inclusion of the \textit{Sigmoid} activation function in some of our networks is essential. However, based on our experiments, we observed that utilizing only the \textit{Sigmoid} activation function in the network can result in slow optimization and convergence. This issue may stem from the problem of vanishing gradient associated with the \textit{Sigmoid} \cite{roodschild2020new,hanin2018neural}.

In response to this challenge, we have thoughtfully incorporated a combination of alternating activation functions into our network architecture. Specifically, after the first layer, which is activated by the \textit{Sigmoid} function, in this work, we set the activation function for odd-numbered layers as \textit{Sigmoid} while setting that in the even-numbered layer as \textit{ReLU}. Besides, it's important to note that the layer before the final output layer does not incorporate any activation function so as to ensure the correct scaling. This deliberate choice of activation functions strikes a balance between ensuring sufficiently continuous derivatives from \textit{Sigmoid}, which are crucial for precise fluid dynamics modeling\cite{jin2021nsfnets}, and the fast convergence offered by the \textit{ReLU} activation\cite{krizhevsky2012imagenet}.

\subsection{Loss Function}

\begin{figure}[!t]
    \centering
    \includegraphics[width=\textwidth]{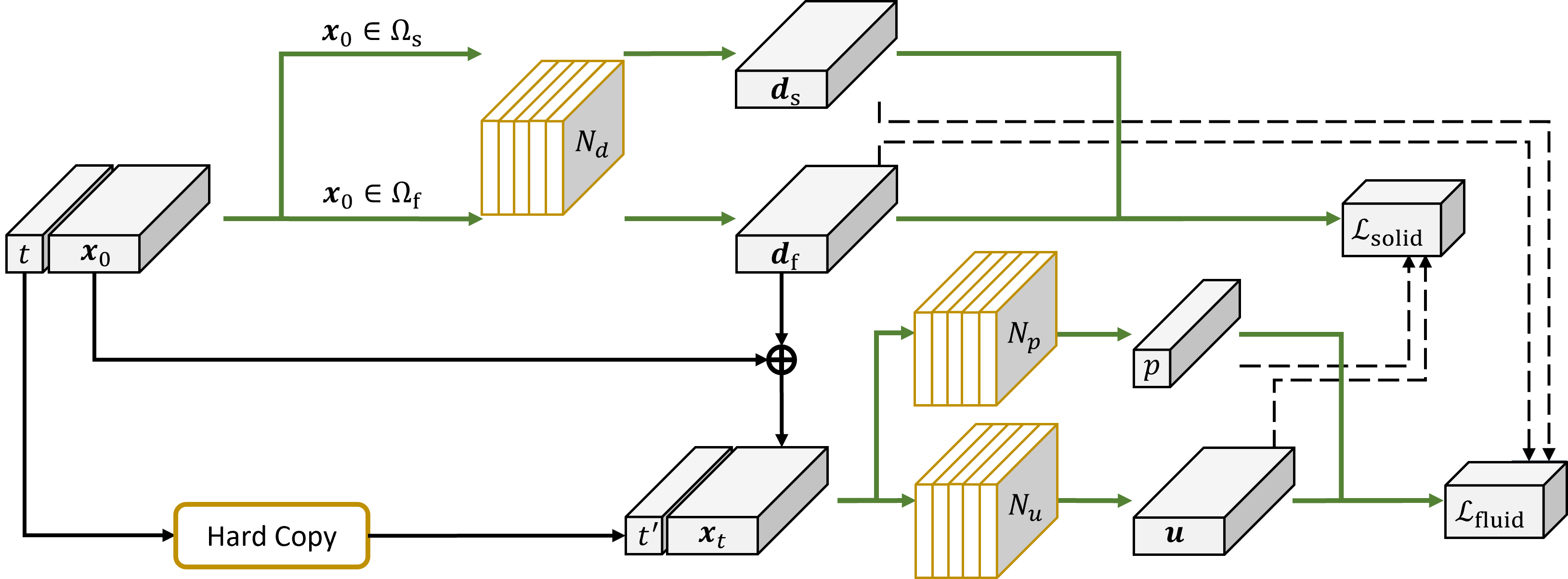}
    \caption{Illustration of Network Architecture: Three networks are employed. The displacement network $N_d$ solves the displacement vector $\boldsymbol{d}$ for both fluid and solid domains. The velocity network $N_u$ computes fluid velocity $\boldsymbol{u}$, while the pressure network $N_p$ handles fluid pressure $p$. Green lines indicate data flow requiring gradient calculation, while black lines denote no gradient delivery. Dotted black lines indicate computed results for loss function calculation without enabling back-propagation, as explained in Section \ref{sec:scheme}.}
    \label{fig:module}
\end{figure}

In our study, we depart from conventional numerical approaches such as Finite Volume and Finite Element methods. Instead, we venture into the intriguing domain of neural networks (NNs) to tackle the complex issue of coupled fluid-structure interaction. At the heart of our methodology lies the utilization of deep neural networks for the approximation of solutions to the Navier-Stokes equations.

Our process initiates with the sampling of points in the reference frame. For each of these points, denoted as $\boldsymbol{x}_0$ at a specific time $t$, we calculate the corresponding solutions for velocity, pressure, and displacement. To achieve this, we introduce a network $N_d$ that takes spatial and temporal coordinates as inputs and outputs the displacement on the corresponding point and time, which is expressed as:
\begin{equation}
    \eta^*(\boldsymbol{x}_0,t; \boldsymbol{\theta}_{d}) = N_d(\boldsymbol{x}_0,t; \boldsymbol{\theta}_{d}) \quad\text{ for }\boldsymbol{x}_0 \in \Omega_t;t\in[0,T].
\end{equation}
and by \eqref{eq:eta}
\begin{equation}
    \boldsymbol{d}^*(\boldsymbol{x}_0,t; \boldsymbol{\theta}_{d}) = \eta^*(\boldsymbol{x}_0,t; \boldsymbol{\theta}_{d}) \boldsymbol{e}_r \quad\text{ for }\boldsymbol{x}_0 \in \Omega_t;t\in[0,T].
    \label{eq:displacement}
\end{equation}

Utilizing the problem coordinates of a point denoted as $\boldsymbol{x}_0$ and the computed displacement vector $\boldsymbol{d}^*(\boldsymbol{x}_0,t)$, we can obtain the corresponding physical coordinates of the point $\boldsymbol{x}_t$ in the current frame at time $t$ as:
\begin{equation}
    \boldsymbol{x}_t^*(\boldsymbol{x}_0,t; \boldsymbol{\theta}_{d}) = \boldsymbol{x}_0 + \boldsymbol{d}^*(\boldsymbol{x}_0,t; \boldsymbol{\theta}_{d}) \quad\text{ for }\boldsymbol{x}_0 \in \Omega_t;t\in[0,T].
    \label{eq:solidnetworkoutput}
\end{equation}

With the acquirement of the physical coordinates, we can get sets of points in the current frame $\Omega_t$. When the problem coordinate of $\boldsymbol{x}_t^*$ is within the fluid domain, we can compute the fluid velocity and pressure on it. To achieve this, we introduce two networks: one for the velocity computation, named $N_u$, and another for  the pressure computation, named $N_p$. The flow velocity and pressure can be calculated as:
\begin{equation}
\begin{aligned}
    \boldsymbol{u}^*(\boldsymbol{x}_t^*,t'; {\boldsymbol{\theta}}_{u}, \boldsymbol{\theta}_{d}) &= N_u(\boldsymbol{x}_t^*,t'; {\boldsymbol{\theta}}_{u})= N_u(\boldsymbol{x}_0 + \boldsymbol{d}^*(\boldsymbol{x}_0,t;\boldsymbol{\theta}_{d}),t';{\boldsymbol{\theta}}_{u})\\
    P^*(\boldsymbol{x}_t^*,t'; {\boldsymbol{\theta}}_{p}, \boldsymbol{\theta}_{d}) &= N_p(\boldsymbol{x}_t^*,t';{\boldsymbol{\theta}}_{p})= N_p(\boldsymbol{x}_0 + \boldsymbol{d}^*(\boldsymbol{x}_0,t;\boldsymbol{\theta}_{d}),t';{\boldsymbol{\theta}}_{p})
\end{aligned}
 \quad\text{ for }\boldsymbol{x}_0 \in \Omega_t^{\text{f}};t\in[0,T].    
 \label{eq:fluidnetworkoutput}
\end{equation}
where for convenience of discussion, we denote $\boldsymbol{\theta}_{f} = \{{\boldsymbol{\theta}}_{u},{\boldsymbol{\theta}}_{p}\}$. 

Please note that in this context, we are using the symbol $t'$, which is distinct from the $t$, to represent a separate hard copy of the value of $t$. Although they have the same value, they are considered different objects in the programming. This distinction is necessary for our implementation using \textit{PyTorch}, as we need to distinguish the differentiation with respect to time in the reference configuration $\Omega_0$ and the current configuration $\Omega_t$.

To clarify, if we were to directly apply automatic differentiation to $\boldsymbol{u}^*(\boldsymbol{x}_t^*;t)$ with respect to $t$, it would yield a result similar to the material derivative as $\frac{\partial \boldsymbol{u}^*}{\partial t} + \frac{\partial \boldsymbol{d}^*}{\partial t} \cdot \nabla \boldsymbol{u}^*$. However, by using the separate variable $t'$, which is a distinct copy of $t$, the automatic differentiation to $\boldsymbol{u}^*(\boldsymbol{x}_t^*;t')$ with respect to $t'$ is $\frac{\partial \boldsymbol{u}^*}{\partial t'}$.

To enable the training of our neural network, we reformulate these partial differential equation (PDE) systems discussed in Section.\ref{sec:problem} into residual forms. In the case of the fluid problem, we rewrite the Navier-Stokes equation and its associated initial/boundary conditions into residual form as:
\begin{equation}
\begin{aligned}
\mathcal{L}_{\text{ns}}({{\boldsymbol{\theta}_{f}}} , {{\boldsymbol{\theta}_{d}}}) &= ||\frac{\partial}{\partial t'} \boldsymbol{u^*}+(\boldsymbol{u}^* \cdot \nabla) \boldsymbol{u}^*+\nabla P^* - 2 \nabla\cdot(\mu \boldsymbol{D}(\boldsymbol{u}^*))||^2_{L^2(\Omega_t^{\text{f}}; [0,T])},\\
&+ ||\nabla\cdot\boldsymbol{u}^*||^2_{L^2(\Omega_t^{\text{f}}; [0,T])}\\
\mathcal{L}_{\text{bdr}}^{\text{f}}({{\boldsymbol{\theta}_{f}}} , {{\boldsymbol{\theta}_{d}}})  &= ||F^{\text{f}}_{\text{bdr}}(\boldsymbol{u}^*,P^*)||^2_{L^2(\Gamma_t^{\text{w}}; [0,T])},\\
\mathcal{L}_{\text{init}}^{\text{f}}({{\boldsymbol{\theta}_{f}}}, \boldsymbol{\theta}_{d})  &= ||F^{\text{f}}_{\text{init}}(\boldsymbol{u}^*,P^*)||^2_{L^2(\Omega_0^{\text{f}}; t=0)}.
\end{aligned}
\label{eqs:fluid}
\end{equation}
Note that the $L^2$ norm is employed here and will be redefined by a discrete approximation in the next subsection. Besides, since we assume the problem to be axis-symmetric, we also include the cylindrical form for the Navier-Stokes equation in Appendix.

Regarding the solid problem, we need to rewrite the stress continuity equation \eqref{eq:StressContinuity}, which is the source of the deformation. The initial and boundary conditions associated are also necessary. The displacement that extends from the displacement on the interface for the fluid domain is also converted for the calculation of the current domain. Those residual forms are:
\begin{equation}
\begin{aligned}
\mathcal{L}_{\text{sc}}({{\boldsymbol{\theta}_{f}}} , {{\boldsymbol{\theta}_{d}}}) &= || \frac{\partial^2 \eta^*}{\partial t^2}+b \eta^* - H||^2_{L^2(\Gamma_0^{\text{w}}; [0,T])},\\
\mathcal{L}_{\text{bdr}}^{\text{s}}({{\boldsymbol{\theta}_{d}}})  &= ||F^{\text{s}}_{\text{bdr}}(\boldsymbol{d}^*)||^2_{L^2( \partial\Gamma_0^{\text{w}}; [0,T])},\\
\mathcal{L}_{\text{init}}^{\text{s}}({{\boldsymbol{\theta}_{d}}})  &= ||F^{\text{s}}_{\text{init}}(\boldsymbol{d}^*)||^2_{L^2(\Gamma_0^{\text{w}}, t=0)},\\
\mathcal{L}_{\text{he}}({{\boldsymbol{\theta}_{d}}}) &= ||\Delta\boldsymbol{d}^*||^2_{L^2(\Omega_0^{\text{f}}; [0,T])};
\end{aligned} 
\label{eqs:solid}
\end{equation}
where $\eta ^*(\boldsymbol{x}_0)$ is the projection of $\boldsymbol{d}^*(\boldsymbol{x}_0)$ onto $\boldsymbol{e}_r$ for $\boldsymbol{x}_0 \in \Gamma_0^{\text{w}}$ as discussed in Sec.\ref{sec:problem}. The cylindrical form for the harmonic extension calculation $\mathcal{L}_{\text{he}}$ is given in Appendix.

Combining all the losses for the flow problem and deformation problem, we get
\begin{align}
&\mathcal{L}_{\text{fluid}}({{\boldsymbol{\theta}_{f}}} , {{\boldsymbol{\theta}_{d}}}) = \alpha_{\text{ns}} \mathcal{L}_{\text{ns}} + \alpha_{\text{bdr}}^{\text{f}} \mathcal{L}_{\text{bdr}}^{\text{f}} + \alpha_{\text{init}}^{\text{f}} \mathcal{L}_{\text{init}}^{\text{f}},\label{eq:fluidloss}\\
&\mathcal{L}_{\text{solid}}({{\boldsymbol{\theta}_{f}}} , {{\boldsymbol{\theta}_{d}}}) = \alpha_{\text{sc}} \mathcal{L}_{\text{sc}} + \alpha_{\text{he}} \mathcal{L}_{\text{he}} + \alpha_{\text{bdr}}^{\text{s}} \mathcal{L}_{\text{bdr}}^{\text{s}} + \alpha_{\text{init}}^{\text{s}} \mathcal{L}_{\text{init}}^{\text{s}}.\label{eq:solidloss}
\end{align}
where $\alpha_{\text{ns}}$, $\alpha_{\text{bdr}}^{\text{f}}$ and $\alpha_{\text{init}}^{\text{f}}$, are corresponding weight parameters for the Navier-Stokes equations, boundary conditions, and initial conditions. In the context of the deformation problem, we utilize $\alpha_{\text{sc}}$, $\alpha_{\text{he}}$, $\alpha_{\text{bdr}}^{\text{s}}$, and $\alpha_{\text{init}}^{\text{s}} $ to assign weights to the loss components associated with the stress continuity, harmonic extension, initial condition, and boundary conditions.

Then, if we can minimize the loss $\mathcal{L}_{\text{fluid}}$ to be approximately zero, the physical systems described in Eq.\eqref{eq:ns} with boundary condition Eq.\eqref{eq:fluidboundary} and initial condition Eq.\eqref{eq:fluidinitial} are approximately satisfied. The solutions provided by networks $N_u$ and $N_p$ have successfully solved the model. Similarly, for $\mathcal{L}_{\text{solid}}$, when it converges sufficiently close to zero, the stress continuity equation in Eq.\eqref{eq:StressContinuity}, along with the appropriate initial and boundary conditions (Eq.\eqref{eq:structureinitial} and Eq.\eqref{eq:structureboundary}), as well as the harmonic extension into the fluid domain, are effectively solved. From our numerical experiments, we found that it is efficient to use the following mutli-objective minmization approach to minimize the two loss functions: 
%our objective neural solvers, namely $N_d$, $N_u$, and $N_p$, are obtained through the following optimization problem:
\begin{equation}
    \min _{{\boldsymbol{\theta}}_d, {\boldsymbol{\theta}_{f}}} \mathcal{L}_{\text{solid}}({{\boldsymbol{\theta}_{f}}} , {{\boldsymbol{\theta}_{d}}}), 
    \label{eq:solidoptimizationformulation}
\end{equation}
and
\begin{equation}
    \min _{{\boldsymbol{\theta}_{f}}, \boldsymbol{\theta}_{d}} \mathcal{L}_{\text{fluid}}({{\boldsymbol{\theta}_{f}}} , {{\boldsymbol{\theta}_{d}}}). 
    \label{eq:fluidoptimizationformulation}
\end{equation}

To decrease the value of the losses, gradient descent methods are applied to update ${\boldsymbol{\theta}_{f}}=\{\boldsymbol{\theta}_{u},\boldsymbol{\theta}_{p}\}$ and ${\boldsymbol{\theta}_{d}}$ using three different optimizer by
\begin{align}
&\boldsymbol{\theta}_{u} \leftarrow \boldsymbol{\theta}_{u} - \tau_{u} \nabla_{\boldsymbol{\theta}_{u}} \mathcal{L}_{\text{fluid}}(\boldsymbol{\theta}_{f} , \boldsymbol{\theta}_{d}),\label{eq:opt_u}\\
&\boldsymbol{\theta}_{p} \leftarrow \boldsymbol{\theta}_{p} - \tau_{p} \nabla_{\boldsymbol{\theta}_{p}} \mathcal{L}_{\text{fluid}}(\boldsymbol{\theta}_{f} , \boldsymbol{\theta}_{d}),\label{eq:opt_p}\\
&\boldsymbol{\theta}_{d} \leftarrow \boldsymbol{\theta}_{d} - \tau_{d} \nabla_{\boldsymbol{\theta}_{d}} \mathcal{L}_{\text{solid}}(\boldsymbol{\theta}_{f} , \boldsymbol{\theta}_{d}),\label{eq:opt_d}
\end{align}
where $\nabla_{\boldsymbol{\theta}_{u}}$, $\nabla_{\boldsymbol{\theta}_{p}}$ and $\nabla_{\boldsymbol{\theta}_{d}}$ are the gradients with respect to $\boldsymbol{\theta}_{u}$, $\boldsymbol{\theta}_{p}$ and $\boldsymbol{\theta}_{d}$. Here, $\tau_{u}$, $\tau_{p}$, $\tau_{d}$ are adaptive step lengths calculated by some methods \cite{amari1993backpropagation,KingBa15}. In the implementation, we utilize three \textit{Adam} optimizers for each of the gradient descents above, whose step lengths are calculated using the method in \cite{KingBa15}.

Our method to solve the above FSI problem is to solve the fluid problem and the solid problem iteratively. In particular, the solution of the fluid problem is obtained by imposing the displacement obtained from solving the solid problem at the interface as boundary conditions. The result for the solid problem is solved by imposing $H$ in the stress continuity \eqref{eq:StressContinuity} on the interface. Correspondingly, in our method, which needs to minimize the two losses, we take an alternative optimization strategy to decouple the FSI problem into two sub-problems. The optimization scheme will be given more specifically in \ref{sec:scheme}.

\subsection{Training Scheme}
\label{sec:scheme}
The fluid-structure interaction problem is particularly challenging due to the coupling of multiple physical phenomena. Adding to this complexity, our approach employs three independent networks: $N_u$, $N_p$, and $N_d$. As a result, devising an effective training strategy for these networks is crucial. Based on the consideration to solve the fluid problem and solid problem alternatively, we present our training scheme in Algorithm \ref{alg:trainingscheme}.

We first apply zero-initialization to the parameters of the output layer of $N_d$, which lets the output all zeros. This initialization allows the solid domain to start in a non-deformed state, aligned with the reference domain and provides a controlled start to the simulation process. Then, with $\alpha_{\text{ns}}=0$, we optimize $N_u$ and $N_p$ for $m_{\text{fluid}}$ epochs alternatively with an alternation every $m_1+m_2$ epochs as indicated by Algorithm \ref{alg:optimizefluid}. This sequence of operations ensures that the networks fulfill the initial and boundary conditions.

After the initialization phase, we further train the networks $N_u$ and $N_p$ to solve the Navier-Stokes equations in the domain of the reference frame. During this stage, we initially set $\alpha_{\text{ns}}=10^{-7}$. Then, we update it by a factor of $10$ when the loss function $\mathcal{L}_{\text{fluid}}$ converges or the training has run for $m_{\text{fluid}}$ epochs. The $\alpha_{\text{ns}}$ update is performed for $E$ times for one simulation problem. This strategy is specifically crafted to address difficulties in decreasing losses when dealing with multiple losses. By weighting the loss terms sequentially, we aim to simplify the optimization process and enhance convergence, allowing the optimizer to focus on more heavily weighted losses.

The final stage of training involves coupling the fluid and solid systems. In this phase, we alternately solve the fluid and solid problems. Specifically, we run $m_{\text{solid}}$ epochs to solve the solid problem for the displacement computation, followed by a maximum of $m_{\text{fluid}}$ epochs for the fluid problem. Either training can be stopped in advance in one alternation if it converges (when the error improvement is less than $\epsilon$ for $100$ epochs). We alternate between these two problems until they converge or run for a maximum of $F$ alternations. The decision to train the networks alternatively is based on the concept of multi-objective optimization, which often employs a variety of alternative techniques to find solutions efficiently. This approach allows for more stable training and better convergence by focusing on optimizing subsets of losses sequentially.

\begin{algorithm}
\caption{Fluid-Structure Interaction Training Scheme}
\label{alg:trainingscheme}
\begin{algorithmic}[1]
\State Initialize $N_d$ with \textit{zero\_initialization};
\State $\alpha_{\text{ns}} \gets 0$;
\State $\boldsymbol{\theta}_{u},\boldsymbol{\theta}_{p} \gets $ Optimize fluid problem using Alg. \ref{alg:optimizefluid} for $m_{\text{fluid}}$ epochs;
\State $\alpha_{\text{ns}} \gets 10^{-8}$;
\For{i=1:5}
    \State $\alpha_{\text{ns}} \gets 10\times\alpha_{\text{ns}}$;
    \State $\boldsymbol{\theta}_{u},\boldsymbol{\theta}_{p} \gets $ Optimize fluid problem using Alg. \ref{alg:optimizefluid} for $m_{\text{fluid}}$ epochs;
\EndFor
\For{i=1:F \textbf{and} $\mathcal{L}_{\text{fluid}}$ doesn't converge}
    \State $\boldsymbol{\theta}_{d} \gets $ Optimize solid problem using \eqref{eq:opt_d} for $m_{\text{solid}}$ epochs;
    \State $\boldsymbol{\theta}_{u},\boldsymbol{\theta}_{p} \gets $ Optimize fluid problem using Alg. \ref{alg:optimizefluid} for $m_{\text{fluid}}$ epochs;
\EndFor
\end{algorithmic}
\end{algorithm}
\begin{algorithm}
\caption{Training Scheme for Fluid Problem}
\label{alg:optimizefluid}
\begin{algorithmic}[1]
\State $k=1$;
\For{$k=1:\frac{m_{\text{fluid}}}{(m1+m2)}$}
    \State $\boldsymbol{\theta}_{u} \gets $ Optimize fluid problem using \eqref{eq:opt_u} for $m_1$ epochs;
    \State $\boldsymbol{\theta}_{p} \gets $ Optimize fluid problem using \eqref{eq:opt_p} for $m_2$ epochs;
\EndFor
\end{algorithmic}
\end{algorithm}

\section{Experiment}
\label{sec:experiment}
This section embarks on a comprehensive evaluation of the efficacy of our proposed methodology through a series of experiments. The primary aim is to demonstrate the robustness and reliability of our approach in a variety of settings. 

Firstly, we undertake a comprehensive assessment of the simulation quality for blood flow within a segment of a cylinder-like vessel. To do a quantitative evaluation, we refer to the calculations by Finite Element Method on a highly dense grid with a small time step as the ground truth. Then, we compare our methods with Finite Element Methods on different grid resolutions and time step sizes according to the relative error. The second facet of our experiments focuses on a self-ablation test to determine the optimal network architecture for fluid problems. We employ networks with varying depths and widths to discern the architecture that best suits our specific setting by simulating the same scenario as the previous blood flow in cylinder vessel simulation. Evaluation of the time for using different numbers of GPUs for parallel training is taken. Also, we expand our evaluation on vessel segments afflicted with a plaque for the third set of experiments. Comparison with Finite Element Method calculations are also performed to provide a holistic view. In the final experiment, we evaluate different sizes of plaque to see how it will influence the stress from the fluid as well as how the blood flow after the stenosis will be affected. 

\subsection{Implementation Details}
\label{sec:implementation}
\subsubsection{Environment and Softwares}
Our methods are implemented across all experiments using the \textit{Python} programming language coupled with the \textit{Pytorch} framework. The Finite Element Method is implemented using the open-source package \textit{FEniCS} and FSI package \textit{turtleFSI}\cite{bergersen2020turtlefsi}. It's noteworthy that the computational meshes are generated using \textit{Gmsh}; however, certain meshes, generated automatically, necessitate manual modifications to make them computable with \textit{FEniCS}. A Linux server equipped with four Tesla A600 Graphics Processing Units is employed to execute these experiments.

\subsubsection{Hyperparameters} 
\label{sec:parameter}
The most important consideration for the weights on the loss functions is ensuring the satisfaction of the boundary and initial conditions. To achieve this, we first train the network to minimize the residuals with $\alpha_{\text{ns}}=0$, $\alpha_{\text{bdr}}^{\text{f}}=1$, and $\alpha_{\text{init}}^{\text{f}}=10^{-1}$ for the fluid problem at the beginning. This setting is motivated by the treatment of initial and boundary conditions in conventional finite element methods, where they are enforced strictly. When the network $N_u$ and $N_p$ are sufficiently trained, we set $\alpha_{\text{ns}}=10^{-7}$ and update it by $\alpha_{\text{ns}} = 10 \times \alpha_{\text{ns}}$ until it converges for the current $\alpha_{\text{ns}}$. This dynamic update helps the neural solver to sufficiently satisfy the initial and boundary conditions before optimizing the network to produce accurate solutions for velocity and pressure. Regarding the solid problem, based on empirical experiments, we set $\alpha_{\text{sc}}=1$, $\alpha_{\text{he}}=10$, $\alpha_{\text{bdr}}^{\text{s}}=10^{-1}$, and $\alpha_{\text{init}}^{\text{s}}=10^{-2}$. 

We set $M=N=1000$ points with temporal and spatial coordinates sampled on each domain and boundary, as discussed in Section \ref{sec:optimization}. Adam optimizer with a learning rate of $1\times 10^{-3}$ is employed to initialize and optimize the learnable variables of neural network. Regarding the network architecture, unless specified otherwise, the network depth is set to $12$ for $N_u$, $N_p$, and $N_d$. The width is set to $20$ for $N_u$ and $N_d$, and $10$ for $N_p$.

For parameters in Algorithm \ref{alg:trainingscheme}, when $\alpha_{\text{ns}}=0$, $N_u$ and $N_p$ are optimized for a maximum of $m_{\text{fluid}}=2000$ epochs which can be stopped when it converges before that. The convergence parameter is $\epsilon= 10^{-1}$. After the initialization stage, we set $\alpha_{\text{ns}}=10^{-8}$, which is then multiplied by $10$ for $E=5$ times. In the phase of iteratively solving the two sub-problems, we run $m_{\text{solid}}=500$ epochs to solve the solid problem, followed by a maximum of $m_{\text{fluid}}=2000$ training units for the fluid problem. There are a maximum of $F=6$ alternations for solving the FSI problem. For parameters in Algorithm \ref{alg:optimizefluid}, training $m_1=80$ epochs for $N_u$ will follow a training on $N_p$ for $m_2=20$ epochs on $N_p$.

It's worth mentioning that the models outlined in Section \ref{sec:problem} are adaptable to dimensional forms with units and non-dimensional forms as in \cite{jin2021nsfnets}. Due to the consideration of applying our methods to real clinical simulations, we retain all units by utilizing the CGS unit system in our implementation. However, it is straightforward to convert the formulation into a non-dimensional form by using the Reynolds number.

\subsubsection{Error Metric}
To evaluate the accuracy and speed of our methods, the relative error and computation time are reported in the sections below. The relative error for velocity magnitude is defined as:
\begin{equation}
    E(\boldsymbol{u}^*,\boldsymbol{u}^{\text{GT}}) = 
    t_l \sum_{\boldsymbol{t}\in\mathcal{T}} \frac{\sum_{\mathcal{P}\in\mathcal{F}}  A(\mathcal{P})
    \frac{\sum_{\boldsymbol{p}\in\mathcal{P}} |\boldsymbol{u}^*(\boldsymbol{p}) - \boldsymbol{u}^{\text{GT}}(\boldsymbol{p})|^2}{|\mathcal{P}|}}{\sum_{\mathcal{P}\in\mathcal{F}}  A(\mathcal{P})\sum_{\boldsymbol{p}\in\mathcal{P}} \frac{|\boldsymbol{u}^{\text{GT}}(\boldsymbol{p})|^2}{|\mathcal{P}|}},
\end{equation}
where $\boldsymbol{u}^*$ represents the solution obtained using our method, while $\boldsymbol{u}^{\text{GT}}$ denotes the ground truth solution, which we utilize the solution obtained using FEM on a high-resolution grid and small step size. The collection of sampled times is denoted by $\mathcal{T}$ with the time step length $t_l$. The set $\mathcal{F}$ represents the collection of elements, while $\mathcal{P}$ denotes one element containing $|\mathcal{P}|$ nodes, and each node is written as $\boldsymbol{p}$. Here $A(\mathcal{P})$ calculates the volume of the element $\mathcal{P}$. For the evaluation of the relative error in pressure, a similar definition is employed.

\subsection{Blood Vessel with Elastic Wall}
\label{sec:problemvessel}
\begin{figure}[ht]
    \centering
    \includegraphics[width=\textwidth]{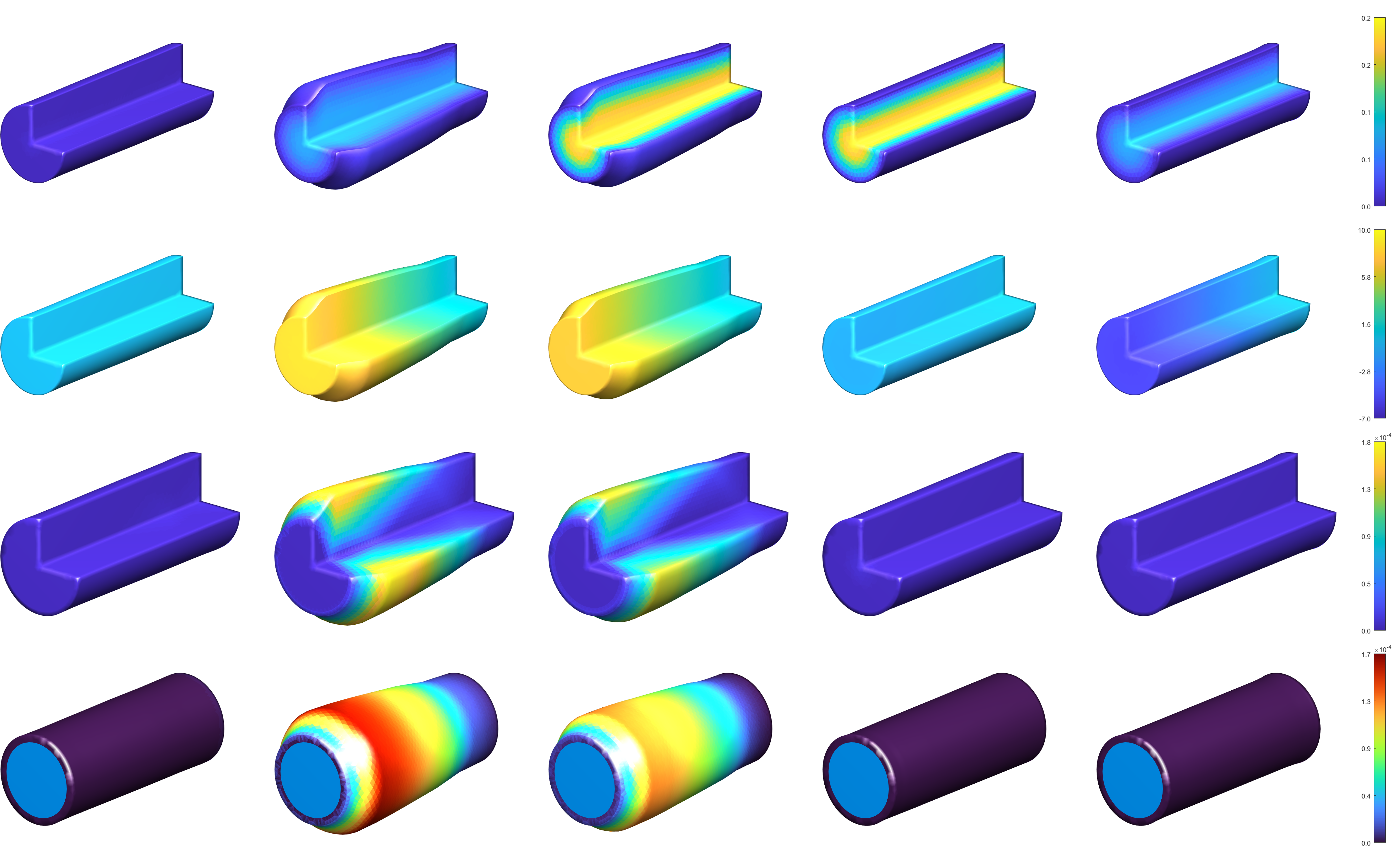}
    \caption{Illustration of velocity magnitude (first row), pressure (second row), displacement magnitude (third row) and displacement magnitude for wall (forth row) variations at time points 0.01, 0.21, 0.41, 0.61, and 0.81. The first and the second don't contain the vessel wall. Deformations of the shape are enlarged by a factor of 10 for better visibility.}
    \label{fig:SequenceVessel}
\end{figure}

In this section, we delve into the details of our simulation process, focusing on the incompressible viscous flow within a linearly elastic vessel. Let us begin by outlining the material parameters and geometry for this study.

The fluid density, denoted by $\rho$, is set to $1.025$ g/$\text{cm}^3$, and the viscosity, denoted as $\mu$, is $0.035$ poise. The vessel has a cylinder-like shape, with a length $l$ of $2$ cm and a radius $r_0$ of $0.25$ cm. The vessel wall thickness is $h=0.05$ cm. The entire simulation runs for a duration $T$ of $2$ seconds.

To elucidate the initial and boundary conditions, let us begin with the initial state of the fluid. At the starting point in time $t=0$, the fluid is at rest, and its velocity $\boldsymbol{u}$ is zero within the fluid domain $\Omega^{\text{f}}_0$, as Equation \eqref{eq:flowinitial}.
\begin{equation}
    \boldsymbol{u}=0 \quad\text{ in }\Omega^{\text{f}}_0, t=0.
    \label{eq:flowinitial}    
\end{equation}

Fluid continuously enters the vessel through the inlet, following a parabolic velocity profile $g(\boldsymbol{x}) = 1-\frac{r^2}{r_0^2}$. Equation \eqref{eq:flowinlet} provides the mathematical representation for this inlet velocity profile. The maximum value of velocity oscillates sinusoidally, achieving a maximum of $20$ cm/s and a minimum of $0$ cm/s, fitting the nature of the cardiac cycle, approximately $60$ times per minute.
\begin{equation}
    \boldsymbol{u}(\boldsymbol{x},t)=g(\boldsymbol{x})(10-10 \cos{(2\pi t)}) \quad\text{ on }\Gamma^{\text{ in }}_t, t\in[0,T].
    \label{eq:flowinlet}    
\end{equation}

On the outlet of the vessel segment, no stress is imposed. This condition is represented as a natural boundary condition, written as
\begin{equation}
    \left[\mu \left( \nabla^\intercal \mathbf{u} + \nabla \mathbf{u} \right)- P \mathbf{I}\right] \cdot \mathbf{n} = \mathbf{0} \quad \text{ on } \ \Gamma_t^{\text{out}}, t \in [0,T].
    \label{eq:flowoutlet}    
\end{equation}

For the interface between the fluid domain and the solid domain, the velocity of the fluid is matched with the speed of the structural movement. Thus, the boundary condition here is as:
\begin{equation}
    \mathbf{u} = \frac{\partial \eta}{\partial t} \mathbf{e}_r \quad \text{ on } \Gamma_t^{\text{w}}, t \in [0,T].
    \label{eq:flowinterface}
\end{equation}

Thus, in the fluid problem, the function $F^{\text{f}}_{\text{init}}$ that indicates the initial condition and $F^{\text{f}}_{\text{bdr}}$ that indicates the boundary condition can be expressed as
\begin{align}
    F^{\text{f}}_{\text{init}} &= \boldsymbol{u} \qquad\qquad\qquad\qquad\qquad\qquad\text{ in }\Omega^{\text{f}}_0, t=0.\\
    \label{eq:functionalinitialfluid}
    F^{\text{f}}_{\text{bdr}} &= 
    \begin{cases}
        \boldsymbol{u}-g(10-10 \cos{(2\pi t)}) \quad&\text{ in }\Gamma^{\text{ in }}_t, t\in[0,T],\\
        \left[\mu \left( \nabla^\intercal \mathbf{u} + \nabla \mathbf{u} \right)- P \mathbf{I}\right] \cdot \mathbf{n} \quad&\text{ on } \Gamma_t^{\text{out}}, t \in [0,T],\\
        \mathbf{u} - \frac{\partial \eta}{\partial t} \mathbf{e}_r \quad&\text{ on } \Gamma_t^{\text{w}},\,\, t \in [0,T].
    \end{cases}
\end{align}

Turning to the parameters for the solid problem, the vessel structure density $\rho_{\text{s}}$ is set at $1.2$ g/$\text{cm}^3$, and its mechanical properties include a Young modulus of $E_{\text{s}}=0.5\times 10^6$ dynes/$\text{cm}^2$ and a Poisson's ratio of $\xi_{\text{s}} = 0.5$ \cite{quarteroni2004mathematical,takashima2007simulation}. Given these mechanical parameters and the governing equation for deformation (Equation \eqref{eq:StressContinuity}), we need suitable initial and boundary conditions to complete the system.

As we discussed in the previous section, the vessel initially maintains a strictly cylindrical shape, ensuring no deformation occurs. This state can be precisely expressed as:
\begin{equation}
    \mathbf{d} = \mathbf{0} \quad\text{ in }\Omega_0, t=0.\label{eq:interfaceinitial}    
\end{equation}

The boundary conditions are established by fixing the starting and ending portions of the vessel segment in their respective positions results in the boundary condition:
\begin{equation}
    \mathbf{d} = \mathbf{0} \quad \text{ on } \partial\Gamma_0^{\text{w}}, t=0.
    \label{eq:interfaceboundary}
\end{equation}

Therefore, the functions $F^{\text{s}}_{\text{init}}$ and $F^{\text{s}}_{\text{bdr}}$, indicating the initial and boundary conditions in the solid problem, can be expressed as:
\begin{align}
    F^{\text{s}}_{\text{init}} &= \boldsymbol{d} \quad\text{ in }\Omega_0, t=0, \label{eq:functionalinitialsolid}\\    
    F^{\text{s}}_{\text{bdr}} &= \mathbf{d} \quad \text{ on } \partial\Gamma_0^{\text{w}}, t=0.     \label{eq:functionalboundarysolid}
\end{align}

\begin{figure}[t!]
    \centering
    % \includegraphics[width=0.4\textwidth]{img/LocationProfileCyn.png}\vspace{3mm}
    % These plots correspond to points indicated in the illustration above, namely, the blue, red, and green points
    \includegraphics[width=\textwidth]{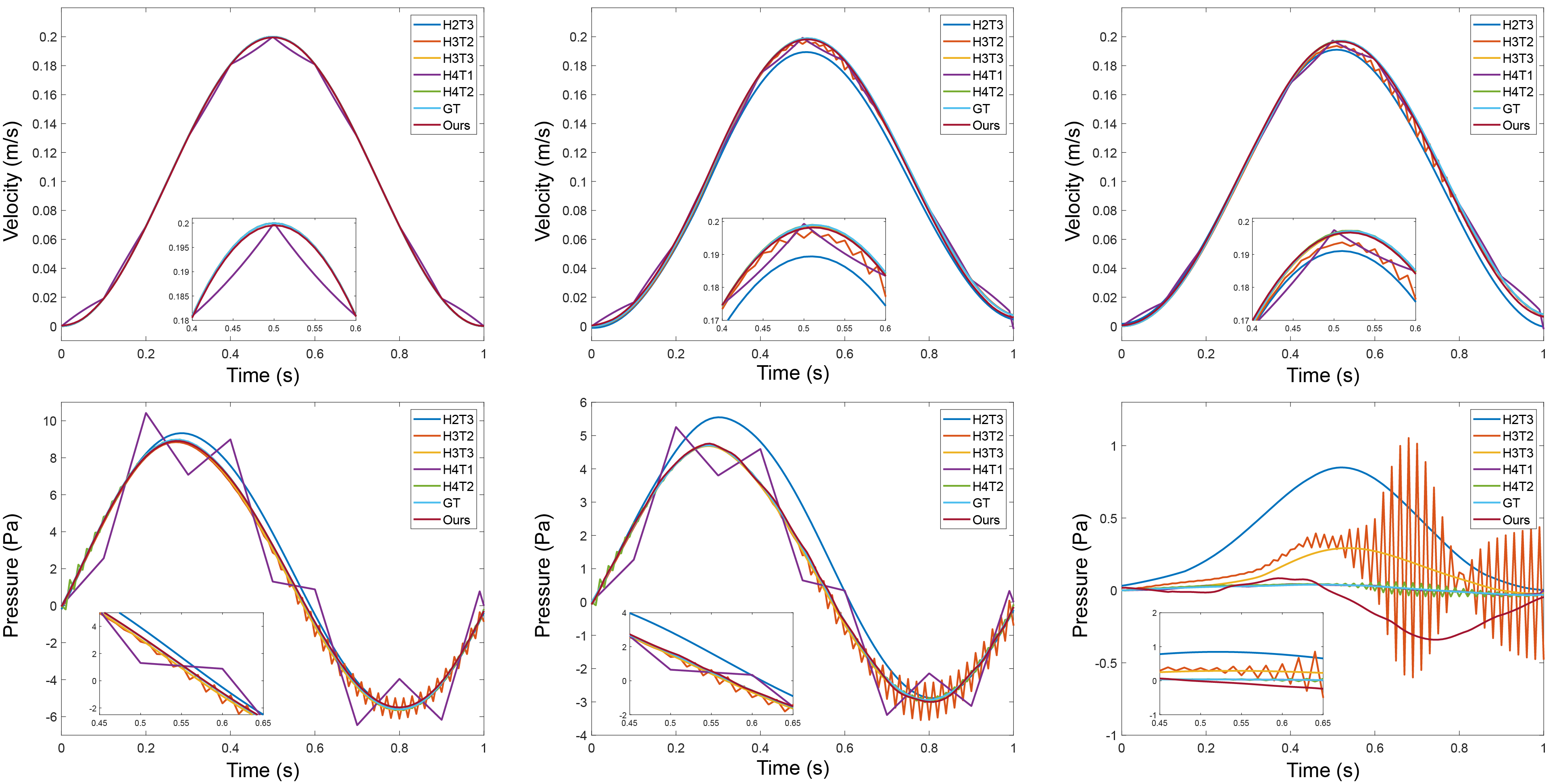}
    \caption{Temporal profiles of velocity magnitude and pressure at different points. From left to right, the upper three plots illustrate velocity magnitude at the centerline's inlet, middle, and outlet, while the lower three represent pressure magnitude for the corresponding points. Legend notations: H1, H2, H3, and H4 indicate grids with minimum lengths of $\frac{1}{2}\times 10^{-3}$, $\frac{1}{4}\times 10^{-3}$, $\frac{1}{8}\times 10^{-3}$, and $\frac{1}{16}\times 10^{-3}$ meter, respectively. T1, T2, and T3 denote time step lengths of $1\times 10^{-1}$, $1\times 10^{-2}$, and $1\times 10^{-3}$ second.}
    \label{fig:VelocityProfileCyn}
\end{figure}
\begin{table*}[t!]
\centering
\begin{tabular}{c|rrrr|r|c}
\toprule
\multicolumn{7}{c}{Relative Error on Velocity Magnitude}\\\hline
$h^{-1}$    &$2\times 10^3$ &$4\times 10^3$ &$8\times 10^3$ & $16\times 10^3$   & Ours  & Points    \\\hline
Time Step $10^{-1}$ & 0.2342    & 0.1089    & 0.0798    & 0.0413    & 0.0363    & $10^2$          \\
Time Step $10^{-2}$ & 0.1142    & 0.0567    & 0.0112    & 0.0008    & 0.0128    & $10^3$          \\
Time Step $10^{-3}$ & 0.0730    & 0.0385    & 0.0003    & 0.0000    & 0.0098    & $10^4$          \\\hline
\hline
\multicolumn{7}{c}{Relative Error on Pressure Magnitude}\\\hline
$h^{-1}$    &$2\times 10^3$ &$4\times 10^3$ &$8\times 10^3$ & $16\times 10^3$   & Ours  & Points    \\\hline
Time Step $10^{-1}$ & 1546.5336 & 4.9901    & 0.0933    & 0.1275    & 0.0415    & $10^2$          \\
Time Step $10^{-2}$ & 17.1869   & 1.4504    & 0.0212    & 0.0017    & 0.0235    & $10^3$          \\
Time Step $10^{-3}$ & 1.9128    & 0.0897    & 0.0083    & 0.0000    & 0.0187    & $10^4$          \\\hline
\hline
\multicolumn{7}{c}{Computation Time (seconds)}\\\hline
$h^{-1}$    &$2\times 10^3$ &$4\times 10^3$ &$8\times 10^3$ & $16\times 10^3$   & Ours  & Points    \\\hline
Time Step $10^{-1}$ & 18        & 51        & 103       & 877       & 732       & $10^2$             \\
Time Step $10^{-2}$ & 87        & 214       & 1053      & 9155      & 863       & $10^3$             \\
Time Step $10^{-3}$ & 283       & 825       & 7168      & 72539     & 2417      & $10^4$            \\\hline
\bottomrule
\end{tabular}
\caption{Comparison between our methods and FEM, reporting relative errors for velocity magnitude and pressure. The FEM result on a high-resolution grid with a small time step is considered ground truth. In the table, $h^{-1}$ signifies the inverse of the minimum length $h$, used to discretize the geometry. Results by our method are shown with different numbers of sample points, corresponding to different time step lengths for FEM.}
\label{tb:comparecos}
\end{table*}

With these critical conditions and material properties defined, we are now fully equipped to proceed with solving the simulation problem for blood flow in elastic vessels, as outlined in Section \ref{sec:method}. The results of our computations are presented in Table \ref{tb:comparecos} where the relative errors for velocity magnitude and pressure are reported. In this comparison, we utilize the solution computed by the Finite Element Method (FEM) on a grid with a minimum length of $\frac{1}{16} \times 10^{-3}$ meter and a time step of $10^{-3}$ seconds as the ground truth. In the table, $h^{-1}$ signifies the inverse of the minimum length $h$, which is used to discretize the geometry. The results for simulation using FEM also relates to time step size, which is set to be three different lengths in our experiments. Notably, our method is mesh-free, so we don't need to configure for time step and discretization length. Instead, we test with a different number of sample points, denoted as $N$, $M$, and $K$ as discussed in Section \ref{sec:method}.

The results indicate that our method performs impressively with $N$, $M$, and $K$ set to $1000$, achieving high accuracy. Remarkably, the accuracy achieved with these settings is comparable to traditional grid-based methods with a minimum length of $1 \times 10^{-4}$ meter and a time step of $10^{-2}$ seconds. However, our approach significantly reduces the computational time.

Figure \ref{fig:SequenceVessel} visually displays the variation in velocity magnitude and pressure at specific time points (0.01, 0.21, 0.41, 0.61, and 0.81). Notably, velocity magnitude is more pronounced near the vessel's centerline, and pressure is generally higher near the inlet, particularly during the systole phase of a cardiac cycle. In contrast, during the diastole phase of a cardiac cycle, the pressure in the front part appears to be negative and lower than at the outlet.

These observations are further supported by the plots in Figure \ref{fig:VelocityProfileCyn}, which depict the temporal profiles of velocity magnitude and pressure at different points. The upper three plots illustrate velocity magnitude, while the lower three represent pressure. To understand the plots, from left to right, they correspond to the points on the centerline's inlet, middle, and outlet. The symbols H1, H2, H3, and H4 in the legend denote the grids whose minimum length $h$ are set to be $\frac{1}{2}\times 10^{-3}$, $\frac{1}{4}\times 10^{-3}$, $\frac{1}{8}\times 10^{-3}$,  and $\frac{1}{16}\times 10^{-3}$ meter, respectively. Similarly, the time step lengths are denoted as T1, T2, and T3. For instance. As an example, the H4T3 corresponds to a grid with a minimum length of $\frac{1}{16}\times 10^{-3}$ meter and a time step of $1 \times 10^{-3}$ second.

Upon closer inspection, we have noticed that the time steps of $1 \times 10^{-1}$ and $1 \times 10^{-2}$ seconds can introduce oscillations that lead to more significant errors over time while our approach consistently produces a smooth solution, which is more accurate and usable in practice.

\begin{figure}
    \centering
    \includegraphics[width=\textwidth]{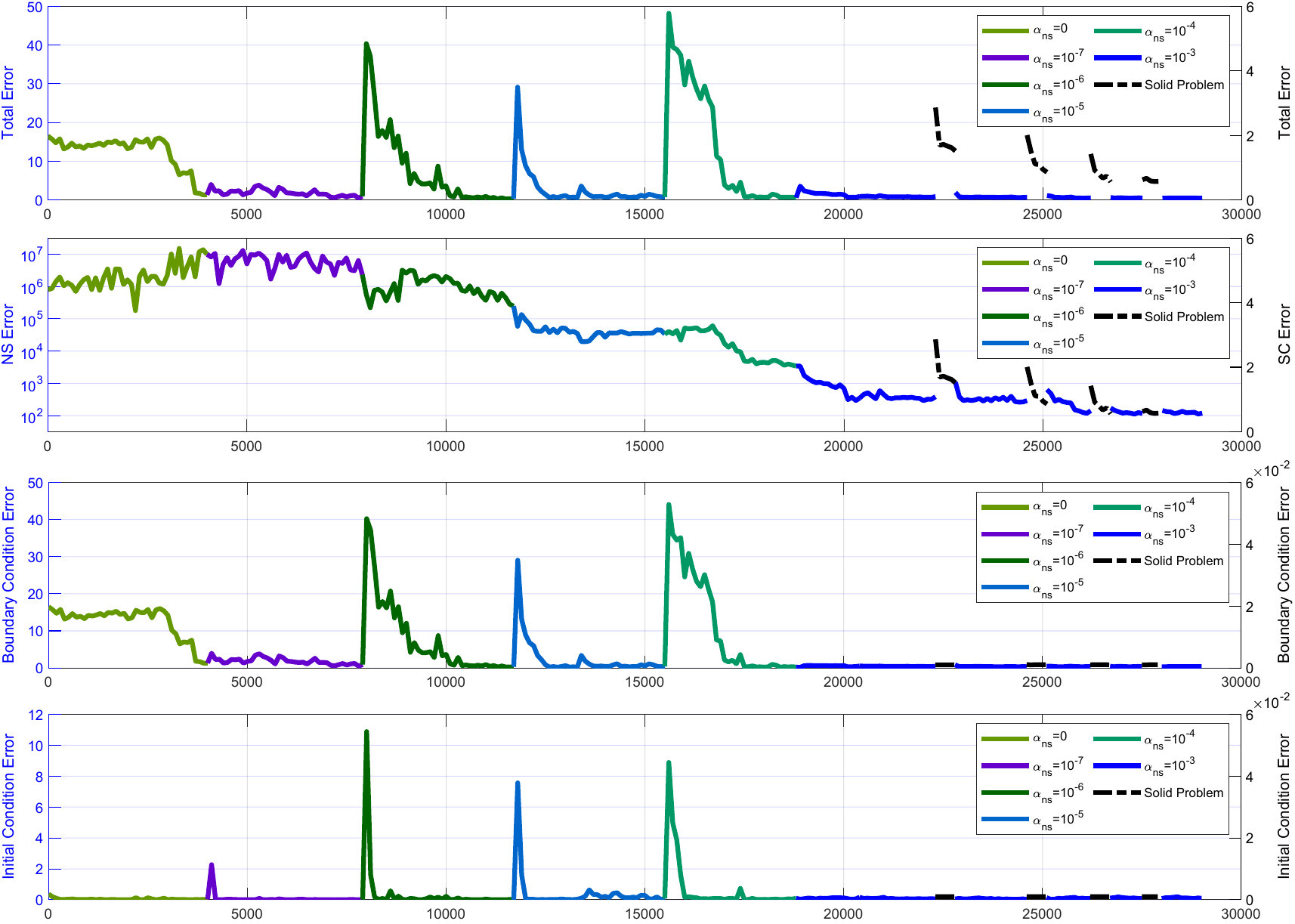}
    \caption{Illustration of loss value variations across training epochs, color-coded to depict changes
in the weighting parameter $\alpha_\text{ns}$ and alternation between fluid and solid problems. Non-black curves represent fluid problem losses at different $\alpha_\text{ns}$ values, with y-axes on the left. Black curves denote solid problem losses, associated y-axes on the right. }
    \label{fig:trainingerror}
\end{figure}

The success of the proposed model greatly benefits from the carefully designed optimization strategy,  which incorporates dynamic weighting parameters and iterative training for both fluid and solid problems. In Figure \ref{fig:trainingerror},we present the variation of the loss values across training for
503 every 100 epochs. The non-black curves represent the loss values for the fluid problem at various $\alpha_\text{ns}$ values, whose y-axes are on the left of the plot. The black curves denote the loss values for the solid problem, whose measure is depicted by the right y-axes.

The left y-axes in blue are the loss values for the fluid problem in different $\alpha_\text{ns}$ values, while the right y-axes in black are the loss values for the solid problem.

In this plot, distinct segments are color-coded to indicate changes in the weighting parameter $\alpha_\text{ns}$. Notably, we observe that as the total loss converges, $\alpha_\text{ns}$ is incrementally updated to a higher value. With each increase in $\alpha_\text{ns}$, the total loss dramatically changed into a high value, prompting the optimizer to subsequently work on reducing the updated error values. For the matching conditions, the associated residual term experiences a substantial surge in value whenever $\alpha_\text{ns}$ is augmented, before being optimized to a lower magnitude. This behavior underscores the efficacy of our dynamic weighting strategy, wherein emphasis is initially placed on improving the matching conditions while concurrently optimizing the Navier-Stokes equation residual.

During the second stage of training, characterized by alternating between the fluid and solid problems (indicated by blue and black segments), the loss steadily diminishes and ultimately converges after four such alternations. Notably, the final training segment for the fluid problem exhibits a marginal reduction in training loss compared to its preceding segment, suggesting a trend toward convergence.

\subsection{Self Ablation}

\begin{figure}[t]
    \centering
    \includegraphics[width=\textwidth]{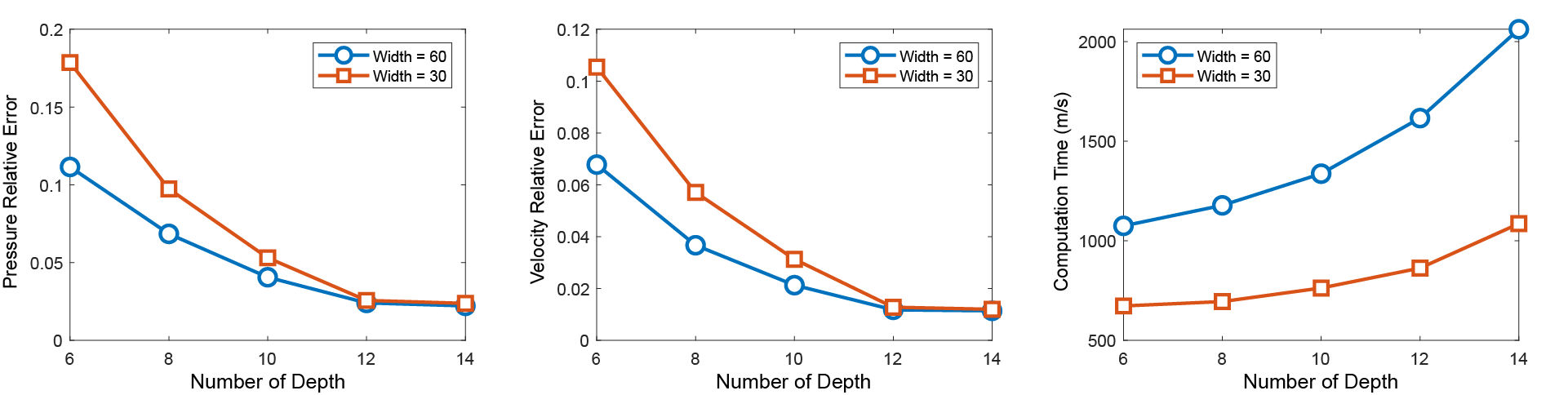}
    \caption{The plots of the solution quality and computational time for networks with different depths and widths. From left to right are pressure relative error, velocity magnitude relative error, and time consumed for reaching a convergence.}
    \label{fig:SelfAblation}
\end{figure}
\begin{table}[ht]
\centering
\small
\begin{tabular}{c|cccccccc}
\toprule
                        &$6\times60$&$8\times60$&$10\times60$&$12\times60$& $14\times60$     \\\hline
No. of Parameters       & 8583      & 12703     & 16823     & 20943     & 25063 \\
Velocity Relative Error & 0.0678    & 0.0367    & 0.0213    & 0.0118    & 0.0114    \\
Pressure Relative Error & 0.1115    & 0.0685    & 0.0406    & 0.0241    & 0.0222    \\
Computational Time      & 1075      & 1178      & 1337      & 1616      & 2062       \\\hline
                        &$6\times30$&$8\times30$&$10\times30$&$12\times30$& $14\times30$      \\\hline
No. of Parameters       & 2293      & 3353      & 4413      & 5473      & 6533\\
Velocity Relative Error & 0.1054    & 0.0571    & 0.0313    & 0.0128    & 0.0120        \\
Pressure Relative Error & 0.1786    & 0.0974    & 0.0531    & 0.0235    & 0.0238        \\
Computational Time      & 673       & 695       & 763       & 863       & 1086 \\\hline
\bottomrule
\end{tabular}
\caption{The velocity magnitude relative error, pressure relative error, and computational time for our methods using different depths and widths.}
\label{tb:selfablation}
\end{table}

\begin{table}[ht]
\centering
\small
\begin{tabular}{c|ccc}
\toprule
                        & Ours      & Sigmoid   & Single NN \\\hline
No. of Parameters       & 5473      & 5473      & 9513      \\
Velocity Relative Error & 0.0128    & 0.0155    & 0.0231    \\
Pressure Relative Error & 0.0239    & 0.0286    & 0.0514    \\
Computational Time      & 863       & 1121      & 1314      \\\hline
\bottomrule
\end{tabular}
\caption{
Comparison of network architectures for fluid problem: Our method uses two FNNs with 12 layers and a width of 30, employing alternating activation functions. ``Sigmoid" refers to an architecture with only Sigmoid activations, while ``Single NN" indicates a single network with alternating \textit{ReLU} and \textit{Sigmoid} activations like ours.}
\label{tb:activation}
\end{table}
\begin{table*}[b!]
\centering
\begin{tabular}{c|cccc}
\toprule
No. of GPUs & $1$   & $2$   & $3$   & $4$ \\\hline
Time (s)    & 863   & 571   & 481   & 440\\
\bottomrule
\end{tabular}
\caption{The computational time for our methods using parallel training on different number of GPUs}
\label{tb:parallel}
\end{table*}

Next, we delve into the influence of different neural network architectures on the model performance. In this section, we maintain the same geometric setup, initial conditions, boundary conditions, and other parameters outlined in the previous simulation problem. The primary variation lies in the neural network architecture.

\textbf{Network Width and Depth} Regarding the depth and the total width for the networks $N_{u}$ and $N_p$, we employ two different width settings. The first setting allocates a width of $20$ to $N_u$ and $10$ to $N_p$, resulting in a total width of $30$. For the second setting, with a total width of $60$, we maintain the same proportion, assigning $40$ to $N_u$ and the remaining $20$ to $N_p$. All tests are conducted under identical optimization settings as detailed in Section \ref{sec:implementation}.

The data presented in Table \ref{tb:selfablation} illustrates that deeper and wider neural networks deliver more precise solutions for velocity and pressure. However, it is worth highlighting that these architectures also impose greater computational consumptions, primarily due to the increased number of trainable parameters. As the networks become deeper and wider, the incremental performance improvement becomes less pronounced, demonstrating a pattern of convergence. Figure \ref{fig:SelfAblation} visually encapsulates these findings. Notably, when the network width is set to be fixed (either 60 or 30), the errors exhibit a discernible convergence pattern. Expanding the number of layers beyond this point does not substantially enhance the relative errors for velocity magnitude and pressure. However, it does significantly extend the duration of the network optimization process.

Considering the trade-off between simulation quality and computational efficiency, we identify the architecture with a depth of 12 and a total width of 30 are the optimal configuration. This architecture is selected for use in the problem of simulating blood flow in a regular cylinder-like vessel segment, as well as the simulation of a vessel with hard plaque within the lumen.

\textbf{Activation Function:} The proposed architecture features an alternating activation function between \textit{ReLU} and \textit{Sigmoid}. While the smoothness of \textit{Sigmoid} is necessary for ensuring at least twice-differentiable solutions, it suffers from vanishing gradient issues away from zero, potentially slowing down convergence. In contrast, \textit{ReLU} lacks this vanishing gradient problem and is computationally more efficient. By alternating between these two functions, we aim to leverage the benefits of both. To assess this approach, we conducted experiments comparing networks with alternating activation functions to those using only \textit{Sigmoid}, employing identical configurations in terms of depth (12 layers) and total width (30 neurons). Results, presented in Table \ref{tb:activation}, demonstrate that the alternating activation function setting yields superior accuracy and significantly reduced computational time.

\textbf{Single vs. Separate Networks:} Existing approaches to solve the Navier-Stokes equation typically employ a single network to predict both pressure and velocity. However, this approach often faced convergence challenges due to the complex internal connections within the hidden layers. Additionally, the Navier-Stokes equation involves distinct operators for velocity and pressure, necessitating tailored approaches for their approximation. As a solution, we employ separate networks for pressure and velocity. This approach provides more flexibility, resulting in improved stability and more precise approximation of the partial differential equations. To validate this, we conducted a comparison experiment on the same simulation problem, contrasting our method with the single-network approach. Results, summarized in Table \ref{tb:activation}, demonstrate that our method achieves superior accuracy in terms of relative errors for both pressure and velocity. Moreover, the optimization time required for convergence is significantly reduced.

\textbf{Parallel Training} Another important feature of the proposed methods is its parallelizability. Benefiting from the well-developed \textit{PyTorch} package, using multiple GPUs to optimize neural networks in a short time is straightforward and simple. Given in Table.\ref{tb:parallel}, we reported the time needed when using different numbers of GPUs for training a network with a depth of $12$ and width of $30$.

\subsection{Blood Vessel with Plaque}
\label{sec:expplaque}
\begin{figure}[t]
    \centering
    \includegraphics[width=0.5\textwidth]{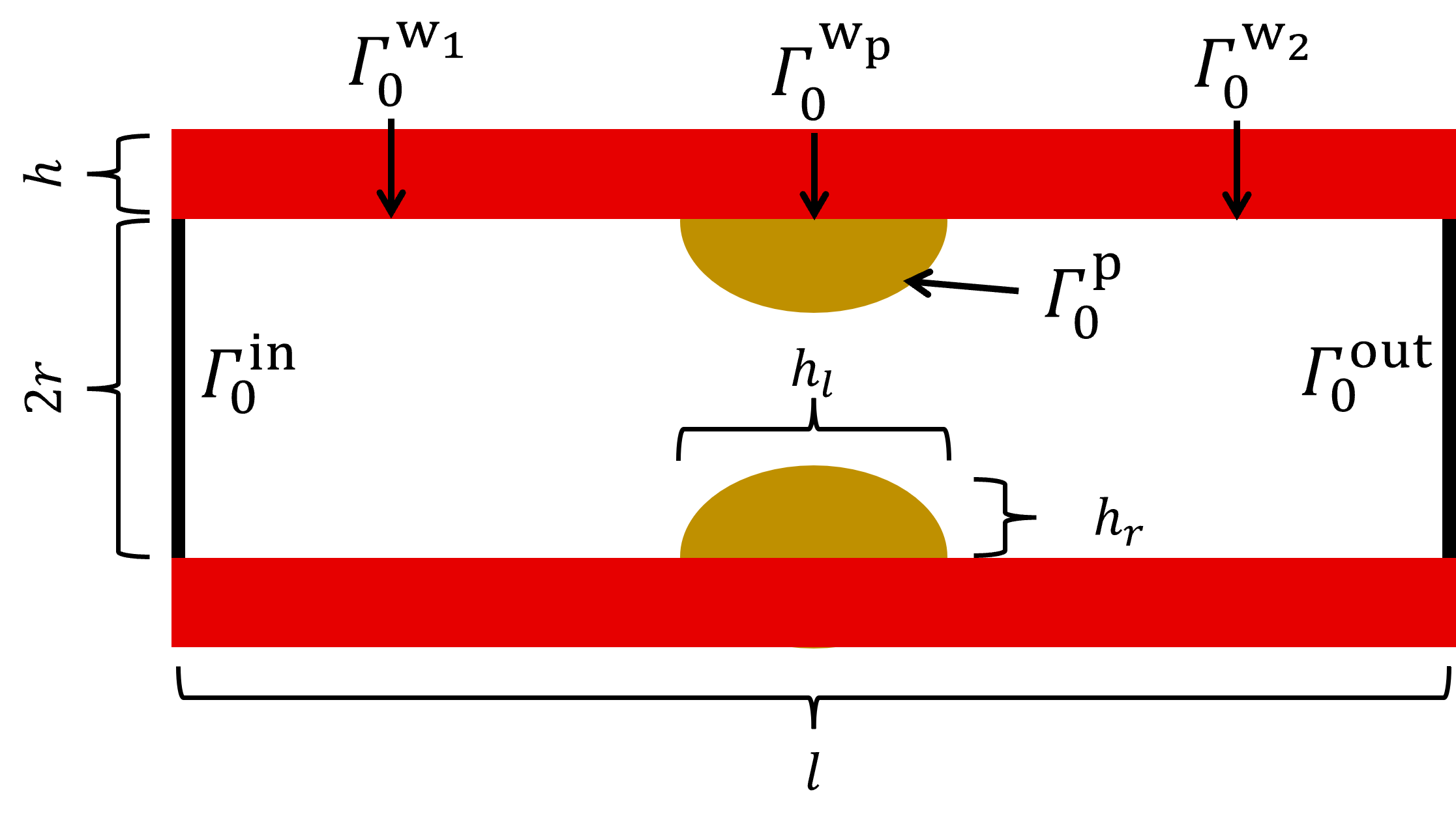}
    \caption{The illustration for the configuration of the cylindrical vessel containing plaque in longitudinal cross-section.}
    \label{fig:GeometryPlaque}
\end{figure}

Here, we are simulating the flow of incompressible viscous fluid within a linearly elastic vessel containing a plaque. The parameters and settings closely resemble those detailed in the first problem discussed in Section \ref{sec:problemvessel}. The parameters of the fluid remain the same, with a density of $\rho_{\text{f}}=1.025$ g/$\text{cm}^3$ and a viscosity $\mu$ of $0.035$ poise.

The vessel retains its cylinder-like shape, as illustrated in Figure \ref{fig:GeometryPlaque}, with a length of $l=2$ cm, a radius of $r_0=0.25$ cm, and a vessel wall thickness of $h=0.05$ cm. In the view of the longitudinal section, the plaque inside the vessel takes the form of a half ellipse centered at the middle point of the vessel wall, as illustrated in Figure \ref{fig:GeometryPlaque}. It has a long radius in the longitudinal direction ($h_l = 0.15$ cm) and a short radius in the radial direction ($h_r = 0.1$ cm). Similar to the previous problem, we simplify this model to a $2$-D setting according to its radial symmetry configuration.

The initial conditions for the fluid remain the same as in the first problem. The fluid starts at rest at the beginning, as described in Equation \eqref{eq:flowinitial}. The function $F^{\text{f}}_{\text{init}}$ from Equation \eqref{eq:functionalinitialfluid} remains applicable in this context. The same inlet and outlet boundary conditions are applied as outlined in Equations \eqref{eq:flowinlet} and \eqref{eq:flowoutlet}. However, the interface is changed into $\Gamma^{\text{w}_1} \cup \Gamma^{\text{w}_2} \cup \Gamma^{\text{p}}$. The continuity of the velocity between the fluid and solid is maintained at the boundaries $\Gamma^{\text{w}_1}_t \cup \Gamma^{\text{w}_2}_t \cup \Gamma^{\text{p}}_t$ by setting the fluid velocity equal to the structure movement speed:
\begin{equation}
    \mathbf{u} = \frac{\partial \eta}{\partial t} \mathbf{e}_r \quad \text{ on } \Gamma^{\text{w}_1}_t \cup \Gamma^{\text{w}_2}_t \cup \Gamma^{\text{p}}_t, t \in [0,T].
\end{equation}

The equation for the function $F^{\text{f}}_{\text{bdr}}$, indicating the boundary condition in the fluid problem, is described as follows:
\begin{equation}
    F^{\text{f}}_{\text{bdr}} = 
    \begin{cases}
        \boldsymbol{u}-g(10-10 \cos{(2\pi t)}) \quad&\text{ in }\Gamma^{\text{ in }}_t, t\in[0,T],\\
        \left[\mu \left( \nabla^\intercal \mathbf{u} + \nabla \mathbf{u} \right)- P \mathbf{I}\right] \cdot \mathbf{n} \quad&\text{ on } \Gamma_t^{\text{out}}, t \in [0,T],\\
        \mathbf{u} - \frac{\partial \eta}{\partial t} \mathbf{e}_r \quad&\text{ on } \Gamma^{\text{w}_1}_t \cup \Gamma^{\text{w}_2}_t \cup \Gamma^{\text{p}}_t, t \in [0,T].
    \end{cases}
\end{equation}

The key difference in this simulation lies in the geometry and mechanical properties introduced by the plaque. Here, we assume the plaque is a hard plaque characterized by the presence of calcium deposits within it\cite{criqui2014calcium}. This type of plaque is rigid and relatively inflexible compared to its softer counterpart. Unlike soft plaques, which can rupture suddenly, hard plaques tend to narrow the arteries, reducing blood flow to the heart muscle and may lead to conditions like coronary artery disease or atherosclerosis.

We retain the parameters for the vessel wall from the previous problem. The vessel structure density is set as $\rho_{\text{s}} = 1.2$ g/$\text{cm}^3$, with Young modulus as $E_{\text{s}} = 0.5\times 10^6$ dynes/$\text{cm}^2$, and a Poisson's ratio of $\xi_{\text{s}} = 0.5$. For the plaque, we set its properties as $\rho_{\text{p}} = 1.1$ g/$\text{cm}^3$ for density, $E_{\text{p}} = 1\times 10^6$ dynes/$\text{cm}^2$ for Young modulus, and a Poisson's ratio of $\xi_{\text{p}} = 0.5$.

The stress continuity equation \eqref{eq:StressContinuity} is modified to account for these different properties:
\begin{equation}
\begin{aligned}
\frac{\partial^2 \eta}{\partial t^2}+b_1 \eta &= H_1 \quad \text{ on } \Gamma^{\text{w}_1}_0 \cup \Gamma^{\text{w}_2}_0,\\
\frac{\partial^2 \eta}{\partial t^2}+b_2 \eta &= H_2 \quad \text{ on } \Gamma^{\text{p}}_0.
\end{aligned}
\label{eq:StressContinuityPlaque}
\end{equation}
where 
\begin{equation}
\begin{aligned}
    b_1 = \frac{E_1}{\rho_{\text{s}_1} (1 - \xi^2 ) R_0^2},
    \quad
    H = \frac{1}{\rho_{\text{s}_1} h_0} \left[ \frac{R}{R_0} P - g \mu \left((\nabla\boldsymbol{u} + (\nabla \boldsymbol{u})^\intercal ) \cdot \boldsymbol{n}\right) \cdot \boldsymbol{e}_r \right],\\
    b_2 = \frac{E_2}{\rho_{\text{s}_2} (1 - \xi^2 ) R_0^2},
    \quad
    H = \frac{1}{\rho_{\text{s}_1} h_0} \left[ \frac{R}{R_0} P - g \mu \left((\nabla\boldsymbol{u} + (\nabla \boldsymbol{u})^\intercal ) \cdot \boldsymbol{n}\right) \cdot \boldsymbol{e}_r \right].
\end{aligned}
\end{equation}
Please note that $R_0$ in \eqref{eq:StressContinuity} has been replaced by $R_p$, a real-value function represented as:
\begin{equation}
R_p(\theta,z) = 
    \begin{cases}
        R_0 \quad &\text{ on } \Gamma^{\text{w}_1}_0 \cup \Gamma^{\text{w}_2}_0,\\
        R_0 - \sqrt{h_r^2-\frac{h_r^2}{h_l^2}(z-z_p)^2} \quad &\text{ on } \Gamma^{\text{p}}_0.
    \end{cases}
\end{equation}

For the initial conditions, the fluid remains at rest initially, consistent with \eqref{eq:interfaceinitial} in the first problem. Consequently, the function $F^{\text{s}}_{\text{init}}$ in \eqref{eq:functionalinitialsolid} is directly applicable in this context. Regarding boundary conditions, as the start and end points of the interface boundary $\Gamma^{\text{w}_1}_0 \cup \Gamma^{\text{p}}_0 \cup \Gamma^{\text{w}_2}_0$ remain the same as those for $\Gamma^{\text{w}}_0$, their descriptions do not differ from what is outlined in \eqref{eq:interfaceboundary}. Thus, we can directly apply the function \eqref{eq:functionalboundarysolid}. 

\begin{figure}[t!]
    \centering
    \includegraphics[width=\textwidth]{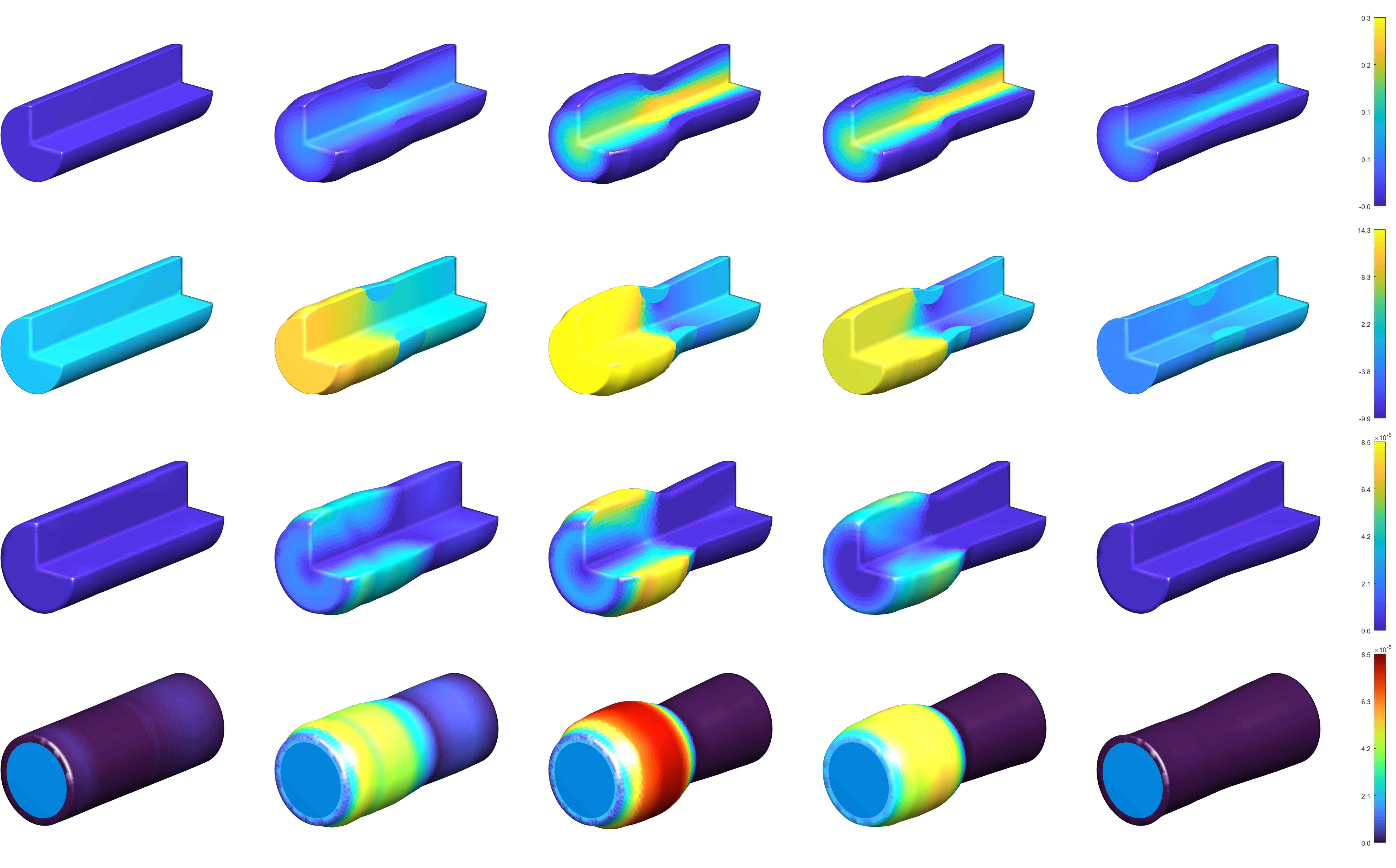}
    \caption{Illustration of velocity magnitude (first row), pressure (second row), displacement magnitude (third row) and displacement magnitude for wall (forth row) variations at time points 0.01, 0.21, 0.41, 0.61, and 0.81. The first and the second don't contain the vessel wall. Deformations of the shape are enlarged by a factor of 10 for better visibility.}
    \label{fig:SequencePlaque}
\end{figure}

We follow the same approach as in the first problem, using a solution obtained from a highly dense grid and a sufficiently small time step as our ground truth. Table \ref{tb:compareplaque} presents the relative errors for velocity magnitude and pressure. As a mesh-free method, we provide results corresponding to different numbers of sample points instead of different time steps and discretization length parameters.

The table shows that our method can deliver accurate results, particularly when using $N$, $M$, and $K$ values greater than $1000$. The accuracy of our method surpasses that achieved by using a grid with a minimum length of $\frac{1}{8} \times 10^{-3}$ meter and a time step of $1 \times 10^{-2}$ seconds. This demonstrates the capability of our method for simulating blood flow in complex geometries.

\begin{figure}[t!]
    \centering
    % \includegraphics[width=0.4\textwidth]{img/LocationProfilePla.png}\vspace{3mm}
    % The plots from left to right correspond to the blue, red, and green points depicted in the accompanying illustration.
    \includegraphics[width=\textwidth]{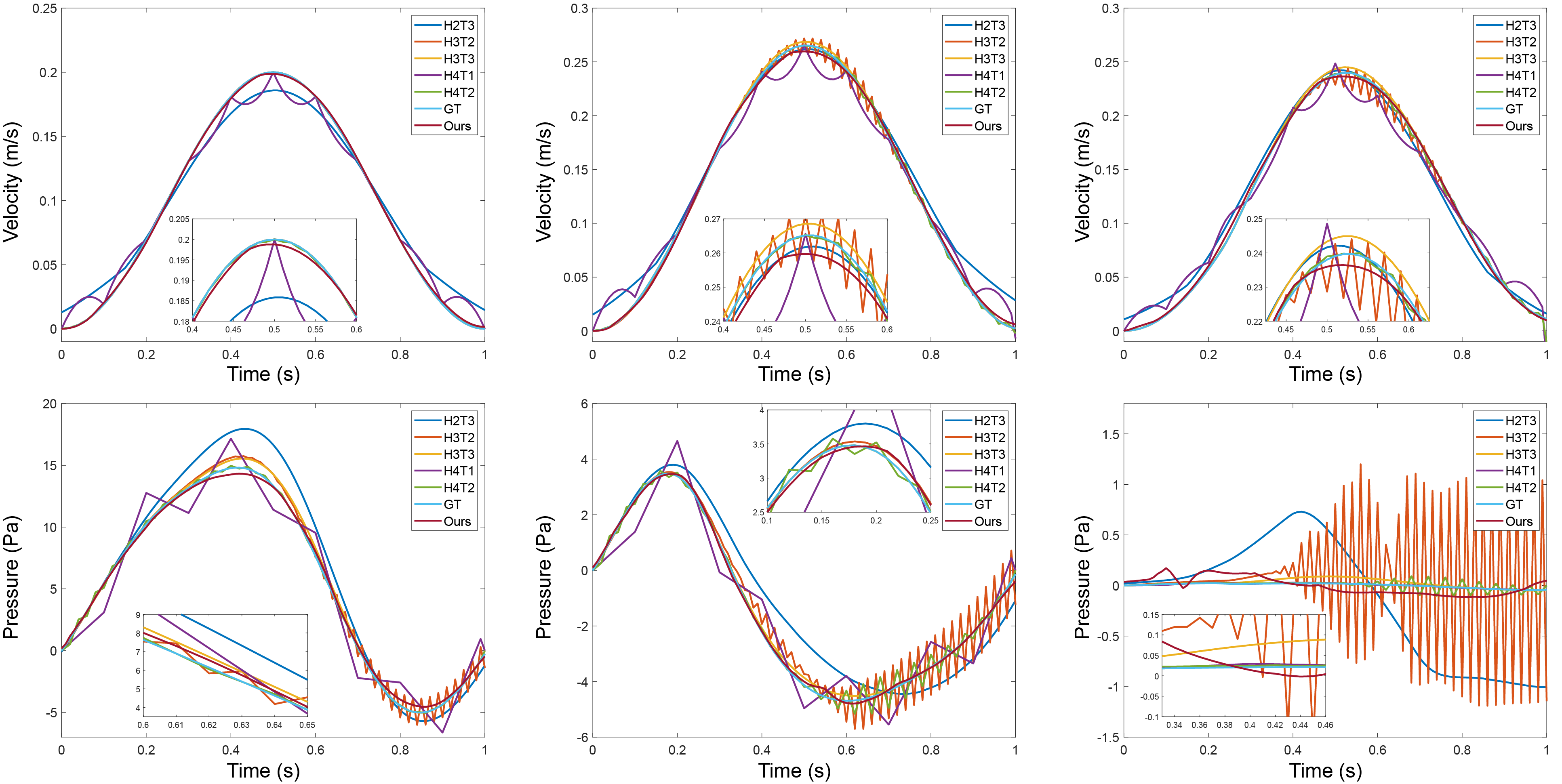}
    \caption{The temporal profiles of the velocity magnitude and pressure at different points. From left to right, the upper three plots illustrate the velocity magnitude on the centerline's inlet, middle, and outlet. The lower three represent the pressure magnitude for the same corresponding points. The symbols H1, H2, H3, and H4 in the legend denote the grids whose minimum length $h$ are set to be $\frac{1}{2}\times 10^{-3}$, $\frac{1}{4}\times 10^{-3}$, $\frac{1}{8}\times 10^{-3}$,  and $\frac{1}{16}\times 10^{-3}$ meter, respectively. Here T1, T2, and T3 denote time step lengths $1\times 10^{-1}$, $1\times 10^{-2}$,  and $1\times 10^{-3}$ second.}
    \label{fig:VelocityProfilePla}
\end{figure}
\begin{table*}[ht]
\centering
\begin{tabular}{c|rrrr|r|c}
\toprule
\multicolumn{6}{c}{Relative Error on Velocity Magnitude}\\\hline
$h^{-1}$    &$2\times 10^3$ &$4\times 10^3$ &$8\times 10^3$ & $16\times 10^3$   & Ours  & Points    \\\hline
Time Step $10^{-1}$ & 2.2252    & 1.7873    & 0.2809    & 0.0515    & 0.0153    & $10^2$    \\
Time Step $10^{-2}$ & 0.2104    & 0.0776    & 0.0136    & 0.0035    & 0.0094    & $10^3$    \\
Time Step $10^{-3}$ & 0.0885    & 0.0390    & 0.0023    & 0         & 0.0092    & $10^4$    \\\hline
\hline
\multicolumn{7}{c}{Relative Error on Pressure Magnitude}\\\hline
$h^{-1}$    &$2\times 10^3$ &$4\times 10^3$ &$8\times 10^3$ & $16\times 10^3$   & Ours  & Points    \\\hline
Time Step $10^{-1}$ & 68.7149   & 41.6116   & 0.9274    & 0.1580    & 0.0378    & $10^2$    \\
Time Step $10^{-2}$ & 22.8401   & 7.3519    & 0.0772    & 0.0153    & 0.0188    & $10^3$    \\
Time Step $10^{-3}$ & 0.8589    & 0.2443    & 0.0073    & 0         & 0.0173    & $10^4$    \\\hline
\multicolumn{7}{c}{Computation Time (seconds)}\\\hline
$h^{-1}$    &$2\times 10^3$ &$4\times 10^3$ &$8\times 10^3$ & $16\times 10^3$   & Ours  & Points    \\\hline
Time Step $10^{-1}$ & 53        & 147       & 411       & 3778      & 2832      & $10^2$             \\
Time Step $10^{-2}$ & 241       & 527       & 4671      & 20512     & 3163      & $10^3$             \\
Time Step $10^{-3}$ & 761       & 3274      & 18375     & 147381    & 9417      & $10^4$            \\\hline
\bottomrule
\end{tabular}
\caption{Comparison between our methods and FEM on cylinder vessel with plaques. The FEM result on a high-resolution grid with a small time step is considered ground truth. $h^{-1}$ signifies the inverse of the minimum length $h$, used to discretize the geometry. Results by our method are shown with different numbers of sample points, corresponding to different time step lengths for FEM.}
\label{tb:compareplaque}
\end{table*}

Following the evaluation in the cylinder-like vessel simulation, we give in Figure \ref{fig:SequencePlaque} the evolution of velocity magnitude and pressure at time points: $0.01$, $0.21$, $0.41$, $0.61$, and $0.81$. The observation reveals a significant increase in velocity magnitude after the stenosis caused by the plaque. This increase results in greater momentum concentration towards the center right between the plaque. This can also be revealed by the plots in Figure.\ref{fig:VelocityProfilePla}, which depict the time histories of the velocity and pressure at the points on the centerline's inlet, middle, and outlet, respectively.

\subsection{Analysis on Plaque with Different Sizes}
Using computational fluid dynamics (CFD) for hemodynamic analysis is an intriguing and highly relevant topic in cardiovascular medicine\cite{li2017numerical}. In this section, we apply our proposed method to investigate the impact of varying plaque sizes. Specifically, we examine four different levels of plaque severity. The configuration for the \textit{No Plaque} case remains the same as the initial problem, featuring a regular, cylinder-like vessel. The other scenarios all involve plaques of different sizes but share similarities with the description provided in Section \ref{sec:expplaque}. To elaborate further, when we focus on the longitudinal cross-section at $\theta = \hat{\theta}$ (as illustrated in Figure \ref{fig:GeometryPlaque}), we set the plaques within the vessel as a half-ellipse centered at the midpoint of the vessel wall. These plaques share a common major radius $h_l = 0.15$ cm in the longitudinal direction and possess varying minor radii $h_r$ in the radial direction, specifically $0.05$ cm for \textit{Mild}, $0.1$ cm for \textit{Moderate}, and $0.15$ cm for \textit{Severe}.

\begin{figure}[h]
    \centering
    \includegraphics[height=4.5cm]{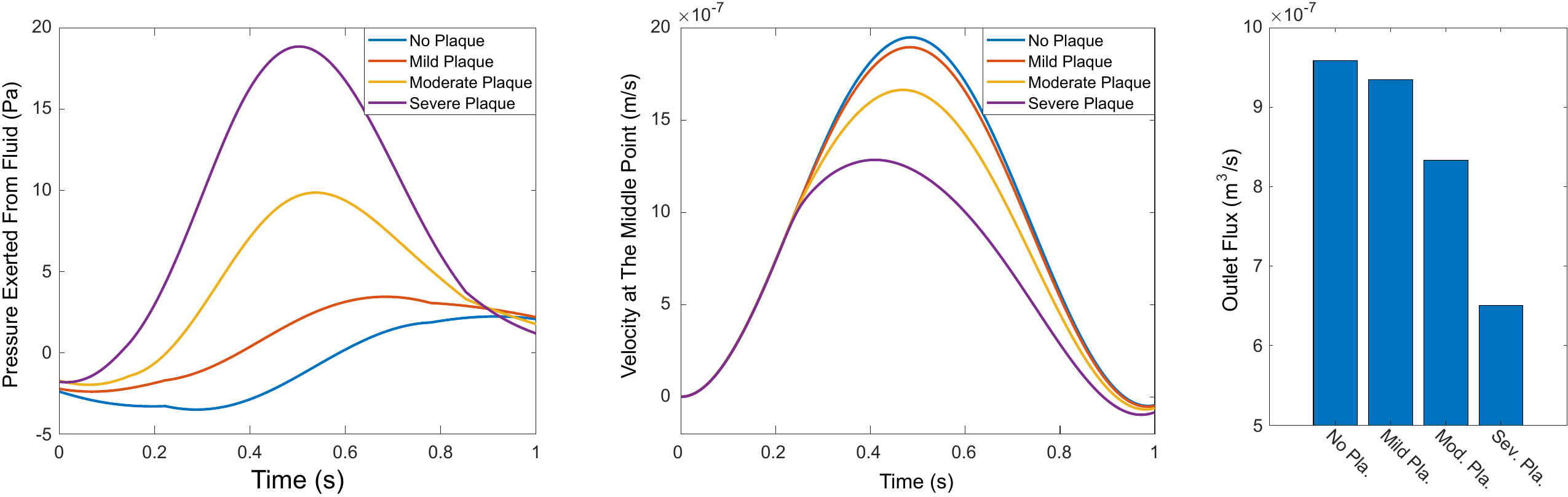}
    \caption{Left: the temporal profiles of the stress on the peak point of the plaque for different plaque sizes; Middle: the temporal profiles of the total flux on the outlet section of the vessel for different plaque sizes; Right: total flux within one cardiac cycle measured on the outlet section of the vessel for different plaque sizes.}
    \label{fig:plaqueflux}
\end{figure}

Firstly, we investigate the effects of force on the vessel for different plaque sizes. In our problem setup, the force on the plaque sourced from the stress exerted from the fluid flow, which is the norm of the traction vector as:
\begin{equation}
    F = ||\boldsymbol{\sigma}_{\text{f}} \cdot \boldsymbol{n}||_2
\end{equation}
where $\boldsymbol{\sigma}_{\text{f}}$ is the Cauchy stress tensor as in \eqref{eq:CompleteStressContinuity}, $\boldsymbol{n}$ is the normal vector to the interface between the plaque and fluid. In the left plot of Figure \ref{fig:plaqueflux}, we present a plot illustrating the force applied by the fluid flow to the highest point of the elliptical plaque. The plot reveals that a healthy vessel without any plaque experiences the least stress in the corresponding region. In contrast, the largest plaque, which significantly alters the vessel geometry, encounters the highest stress from the blood flow. As the minor radius $h_r$ of the elliptical plaque increases linearly, the maximum stress within one cardiac cycle rises exponentially. This demonstrates that an increase in the amount of deposited materials within the plaque leads to a much more rapid and significant escalation in the level of risk.

Next, we explore the impact of plaque on blood flow, specifically on how different plaque sizes affect the flow rate. As shown in Figure \ref{fig:plaqueflux}, a healthy vessel without plaque allows ample blood flow to pass through the outlet boundary of the vessel segment. However, vessels with plaques impede the flow, with larger plaques causing more pronounced restrictions. Similarly, the increase in the short radius of the elliptical plaque results in a notably more dramatic reduction in blood flow.

\section{Conclusion}
In this paper, we proposed a mesh-free approach to simulate the blood flow in a deformable vessel. The blood flow is formulated as the Navier-Stokes equation associated with proper initial and boundary conditions. The mechanical model for the vessel wall structure is modeled by equating Newton's second law of motion and linear elasticity to the force exerted by the fluid flow, which we referred to as the stress continuity equation \eqref{eq:CompleteStressContinuity}. Our method utilizes the physics-informed neural network (PINN) to solve the Navier-Stokes equation coupled with stress continuity on the interface of vessel structure and fluid flow. As a mesh-free approach, our method does not necessitate discretization and meshing of the computational domain and has been proven highly efficient in solving simulations for complex geometries. With the well-developed neural network packages and parallel modules, our method can be easily accelerated through GPU computing and parallel computing. Furthermore, our model can easily be adapted to arteries in different geometries without requiring additional discretization. Therefore, our method is significantly more efficient and user-friendly than employing the conventional Finite Element Method for clinical purposes. To evaluate our methods, we performed experiments on vessels with and without plaque inside. The results indicate the advantages of our methods. Self-ablation is also performed to determine the optimal setting for network architecture. The time saved by parallel training is evaluated. Additionally, we analyze the influence of different plaque sizes using our method. 

In our future work, we aim to leverage Physics-Informed Neural Networks (PINNs) for the inverse design of optimal stent shapes in patient-specific arteries, aiding in the clinical treatment of arterial stenosis \cite{lu2021physics}. The interaction to non-linear elastic and porous solid structure are also important to model. Additionally, given the current capabilities of imaging techniques to capture precise and high-resolution cardiovascular vessel deformations, investigating flow simulation within dynamic domains presents an intriguing and practical avenue for research.

When considering network architectures for our PINN models, we plan to explore more advanced and innovative designs. For instance, integrating principles of domain decomposition, as demonstrated in approaches like xPINN \cite{shukla2021parallel,jagtap2020extended}, could enhance the accuracy of our results, particularly when addressing challenges beyond the capabilities of other methods. Furthermore, exploring architectures that employ adaptive activation functions \cite{jagtap2020adaptive} instead of the alternating \textit{Sigmoid} and \textit{ReLU} setup could potentially lead to improved convergence rates and overall performance.
\bibliographystyle{siamplain} 
\bibliography{references}
\newpage

\appendix
\section{One-Pulse Flow Simulation}
\begin{figure}[ht!]
    \centering
    \includegraphics[width=\textwidth]{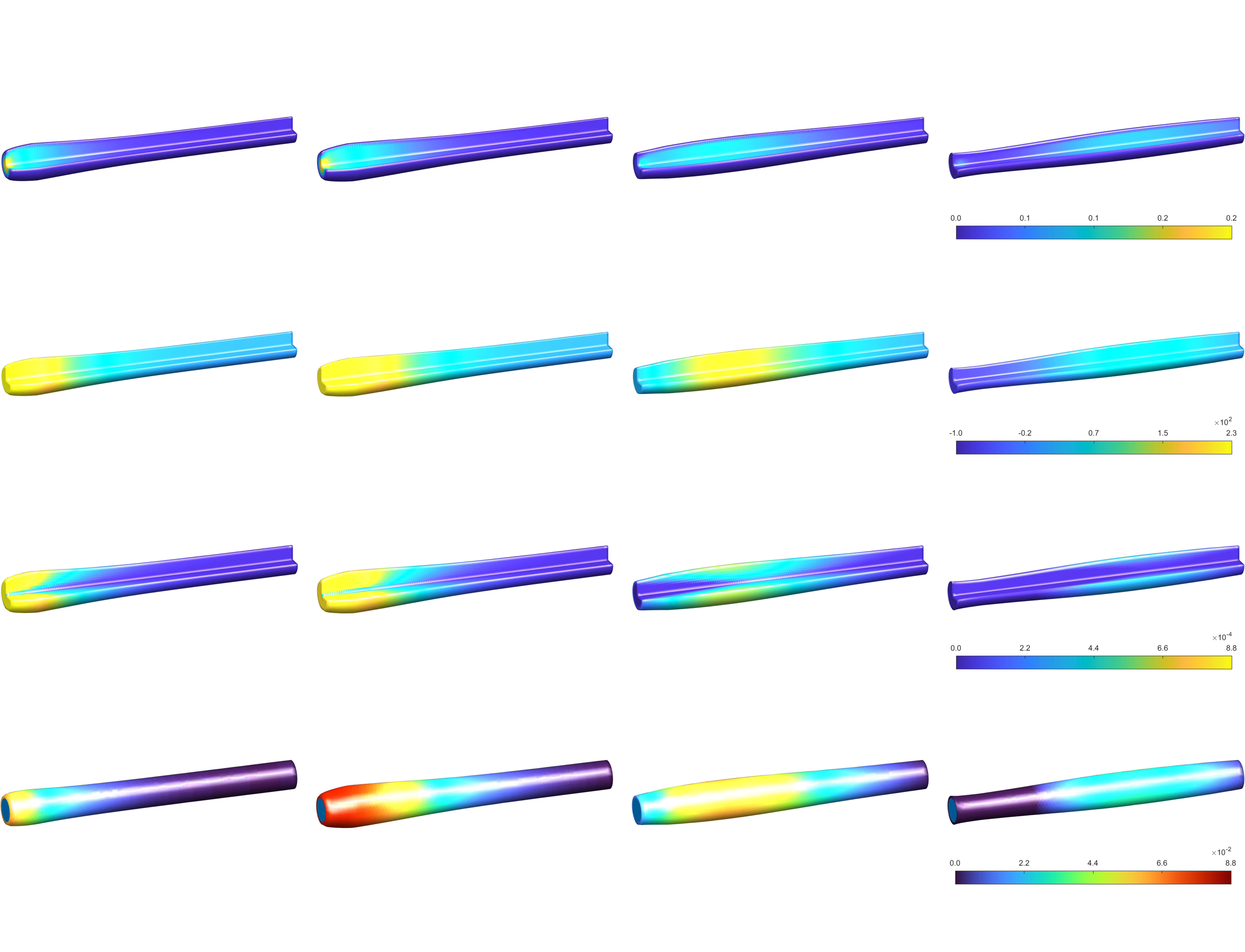}
    \caption{Illustration of velocity magnitude (first row), pressure (second row), displacement magnitude (third row) and displacement magnitude for wall (forth row) variations at time points 0.04, 0.07, 0.10, and 0.14. The first and the second don't contain the vessel wall. Deformations of the shape are enlarged by a factor of 10 for better visibility.}
    \label{fig:propogate}
\end{figure}

Here we perform an experiment with an impulse input at the inlet of the vessel as the inlet boundary condition so that a dilation phenomenon that propagates downward to the outlet can observed. Within this experiment, the fluid density $\rho$ is set to $1.025$ g/$\text{cm}^3$. The viscosity $\mu$ is $0.035$ poise. The vessel has a cylinder-like shape, with a length $l$ of $25$ cm and a radius $r_0$ of $1$ cm. The vessel wall thickness is $h=0.2$ cm. The vessel structure density is set as $\rho_{\text{s}} = 1.2$ g/$\text{cm}^3$, with Young modulus as $E_{\text{s}} = 0.8\times 10^7$ dynes/$\text{cm}^2$, and a Poisson's ratio of $\xi_{\text{s}} = 0.5$. The entire simulation runs for a duration $T$ of $0.2$ seconds.

Fluid enters the vessel through the inlet, following a parabolic velocity profile $g(\boldsymbol{x}) = 1-\frac{r^2}{r_0^2}$. Equation \eqref{eq:flowinlet_propogate} provides the mathematical representation for this inlet velocity profile. The impulse input increases and decreases very quickly so that the flow will build up a segment of high pressure part that flows down. In our setting, the velocity speed at the input will finish increasing and decreasing within 0.1 second, which is described as:
\begin{equation}
    \boldsymbol{u}(\boldsymbol{x},t)=g(\boldsymbol{x})(10-10 \cos{(20\pi t)}) \quad\text{ on }\Gamma^{\text{ in }}_t, t\in[0,T].
    \label{eq:flowinlet_propogate}    
\end{equation}

Other settings like outlet boundary conditions are kept as same as those in the previous cylinder vessel simulation problem. The solved result is visualized in Figure.\ref{fig:propogate}, whose layout and color visualized follow the same rule as those in the main context. The frames rendered from left to right are at times 0.04, 0.07, 0.10, and 0.14. The dilation propagation can be clearly observed.

\section{Navier-Stokes Equation in Cylindrical System}
In our problem, since the fluid is incompressible, we can have $\boldsymbol{\sigma}_f = P \boldsymbol{ I} + 2\mu \nabla \boldsymbol{u}$ as $\nabla\cdot (\nabla\boldsymbol{u})^T=0$ according to \cite[P.31, Sec. 8]{quarteroni2004mathematical}. Then, the Navier-Stokes equation in polar coordinates takes the following form
\begin{equation}
\begin{aligned}
\rho_\text{f}\frac{\partial u_z}{\partial t} + \rho_\text{f}(u_r \frac{\partial}{\partial r} + \frac{u_\theta}{r} \frac{\partial}{\partial \theta} + u_z \frac{\partial}{\partial z}) u_z &= - \frac{\partial P}{\partial z} + \mu \nabla^2 u_z \\
\rho_\text{f}\frac{\partial u_r}{\partial t} + \rho_\text{f}(u_r \frac{\partial}{\partial r} + \frac{u_\theta}{r} \frac{\partial}{\partial \theta} + u_z \frac{\partial}{\partial z}) u_r -\frac{u_\theta^2}{r} &= -\frac{\partial P}{\partial r} + \mu \left( \nabla^2 u_r - \frac{u_r}{r^2} - \frac{2}{r^2} \frac{\partial u_\theta}{\partial \theta} \right)\\
\rho_\text{f}\frac{\partial u_\theta}{\partial t} + \rho_\text{f}(u_r \frac{\partial}{\partial r} + \frac{u_\theta}{r} \frac{\partial}{\partial \theta} + u_z \frac{\partial}{\partial z}) u_\theta + \frac{u_r u_\theta}{r} &= -\frac{1}{r} \frac{\partial P}{\partial \theta} + \mu \left( \nabla^2 u_\theta + \frac{u_\theta}{r^2} + \frac{2}{r^2} \frac{\partial u_r}{\partial \theta} \right)   
\end{aligned}
\label{eq:B1}
\end{equation}
where $\nabla^2 = \frac{1}{r} \frac{\partial}{\partial r}\left(r \frac{\partial }{\partial r}\right)+\frac{1}{r^2} \frac{\partial^2 }{\partial \theta^2}+\frac{\partial^2 }{\partial z^2}$. The incompressibility reads,
\begin{equation}
    \frac{1}{r}\frac{\partial}{\partial r} (r u_r) + \frac{1}{r} \frac{\partial u_\theta}{\partial \theta} + \frac{\partial u_z}{\partial z} = 0
     \label{eq:B2}
\end{equation}
and the harmonic extension is as:
\begin{equation}
\begin{aligned}
    \frac{1}{r} \frac{\partial}{\partial r} (r \frac{\partial u_z}{\partial r} ) + \frac{1}{r^2} \frac{\partial^2 u_z}{\partial \theta^2} + \frac{\partial^2 u_z}{\partial z^2} = 0 \\
    \frac{1}{r} \frac{\partial}{\partial r} (r \frac{\partial u_r}{\partial r} ) + \frac{1}{r^2} \frac{\partial^2 u_r}{\partial \theta^2} + \frac{\partial^2 u_r}{\partial z^2} = 0 \\
    \frac{1}{r} \frac{\partial}{\partial r} (r \frac{\partial u_\theta}{\partial r} ) + \frac{1}{r^2} \frac{\partial^2 u_\theta}{\partial \theta^2} + \frac{\partial^2 u_\theta}{\partial z^2} = 0 
\end{aligned}
 \label{eq:B3}
\end{equation}

As we assume that our problem is axis-symmetric, we have $u_\theta = 0$ for all $\theta$. Accordingly, \eqref{eq:B1} is reduced into
% \begin{equation}
% \begin{aligned}
% \rho_\text{f}\frac{\partial u_r}{\partial t} + \rho_\text{f}(u_r \frac{\partial}{\partial r} + u_z \frac{\partial}{\partial z}) u_r &= - \frac{\partial P}{\partial r} + \mu \left( \frac{1}{r} \frac{\partial}{\partial r} (r \frac{\partial u_r}{\partial r} ) + \frac{\partial^2 u_r}{\partial z^2} - \frac{u_r}{r^2} \right)\\
% \rho_\text{f}\frac{\partial u_z}{\partial t} + \rho_\text{f}(u_r \frac{\partial}{\partial r} + u_z \frac{\partial}{\partial z}) u_z &= - \frac{\partial P}{\partial z} + \mu \left( \frac{1}{r} \frac{\partial}{\partial r} (r \frac{\partial u_z}{\partial r} ) + \frac{\partial^2 u_z}{\partial z^2} \right)
% \end{aligned}
% \end{equation}
\begin{equation}
\begin{aligned}
\rho_\text{f}\frac{\partial u_z}{\partial t} + \rho_\text{f}(u_r \frac{\partial}{\partial r} + u_z \frac{\partial}{\partial z}) u_z &= - \frac{\partial P}{\partial z} + \mu \left( \frac{1}{r} \frac{\partial u_z}{\partial r} + \frac{\partial^2 u_z}{\partial r ^2} + \frac{\partial^2 u_z}{\partial z^2} \right)\\
\rho_\text{f}\frac{\partial u_r}{\partial t} + \rho_\text{f}(u_r \frac{\partial}{\partial r} + u_z \frac{\partial}{\partial z}) u_r &= - \frac{\partial P}{\partial r} + \mu \left( \frac{1}{r} \frac{\partial u_r}{\partial r} + \frac{\partial^2 u_r}{\partial r ^2} + \frac{\partial^2 u_r}{\partial z^2} - \frac{u_r}{r^2} \right)
\label{eq:nscylinderical}
\end{aligned}
\end{equation}
and the incompressibility \eqref{eq:B2} is reduced into 
\begin{equation}
    \frac{1}{r} u_r + \frac{\partial u_r}{\partial r} + \frac{\partial u_z}{\partial z} = 0
    \label{eq:ns2cylinderical}
\end{equation}
and the harmonic extension \eqref{eq:B3} is reduced into 
\begin{equation}
\begin{aligned}
    \frac{1}{r} \frac{\partial u_z}{\partial r} + \frac{\partial^2 u_z}{\partial r ^2} + \frac{\partial^2 u_z}{\partial z^2} = 0 \\
    \frac{1}{r} \frac{\partial u_r}{\partial r} + \frac{\partial^2 u_r}{\partial r ^2} + \frac{\partial^2 u_r}{\partial z^2} = 0 \\ 
\end{aligned}
\end{equation}

The cylindrical form of a PDE would introduce $\frac{1}{r}$ which results in singularities at $r=0$. To deal with it, we first map $[0,r_0]\times[0,2\pi)$ into $[-r_0,r_0]\times[0,\pi)$ and let $r' = sgn(r)(ReLU(|r|-\epsilon_r)+\epsilon_r)$ so that $|\frac{1}{r'}|$ is bounde by $\frac{1}{\epsilon_r}$ near $r=0$, where $sgn(r)$ denotes the sign of $r$. Then the residual form in the cylindrical representation for Navier-Stokes equation is as

\begin{equation}
\begin{aligned}
    &\mathcal{L}_{\text{ns}}({{\boldsymbol{\theta}_{f}}} , {{\boldsymbol{\theta}_{d}}}) \\
    &= ||\rho_\text{f}\frac{\partial u_z^*}{\partial t} + \rho_\text{f}(u_r^* \frac{\partial}{\partial r} + u_z^* \frac{\partial}{\partial z}) u_z^* + \frac{\partial P^*}{\partial z} - \mu \left( \frac{1}{r'} \frac{\partial u_z^*}{\partial r} + \frac{\partial^2 u_z^*}{\partial r ^2} + \frac{\partial^2 u_z^*}{\partial z^2} \right)||^2_{L^2(\Omega_t^{\text{f}}; [0,T])}\\
    &+ ||\rho_\text{f}\frac{\partial u_r^*}{\partial t} + \rho_\text{f}(u_r^* \frac{\partial}{\partial r} + u_z^* \frac{\partial}{\partial z}) u_r^* + \frac{\partial P^*}{\partial r} - \mu \left( \frac{1}{r'} \frac{\partial u_r^*}{\partial r} + \frac{\partial^2 u_r^*}{\partial r ^2} + \frac{\partial^2 u_r^*}{\partial z^2} - \frac{u_r^*}{r'^2} \right)||^2_{L^2(\Omega_t^{\text{f}}; [0,T])}\\
    &+ ||\frac{1}{r'} u_r^* + \frac{\partial u_r^*}{\partial r} + \frac{\partial u_z^*}{\partial z}||^2_{L^2(\Omega_t^{\text{f}}; [0,T])}
\end{aligned}
\end{equation}
where $(u_z^*,u_r^*) = \boldsymbol{u}^*$ and $P^*$ are calculated from \eqref{eq:fluidnetworkoutput}, whose input $\boldsymbol{x}_t^*$ is calculated from \eqref{eq:solidnetworkoutput} by the input $\boldsymbol{x}_0= (r,z)$. We set $\epsilon_r = \frac{r_0}{100}$ in the implementation.

\section{Mini-batch gradient descent algorithm}\label{sec:discrete}
\label{sec:optimization}
Our FSI problem contains the fluid problem \eqref{eq:fluidloss} and the solid problem \eqref{eq:solidloss}. Each problem contains multiple losses defined as the norm in individual domains. To obtain those norms, we follow the principle of \textit{PINN} to sample a sufficient number of points in the domain and treat them as an approximate representation of it. For example, we sample a set of $M$ points in domain $\Omega_0^{\text{f}}$ to obtain $\mathcal{S}({\Omega_0^{\text{f}}})$, which is an approximate representation of it. We additionally sample $M$ temporal points $\mathcal{T}([0,T])$ in $[0,T]$ as an approximate representation of this time interval. Then, we can have 
\begin{equation}
\begin{aligned}
    \mathcal{S}({\Omega_0^{\text{f}}};t=0) = \{ (\boldsymbol{x}^{(i)}_{0},0) \quad|\quad \text{ where }  \boldsymbol{x}^{(i)}_{0}\in\mathcal{S}({\Omega_0^{\text{f}}})\}    
\end{aligned}
\label{eqnorm:t0}
\end{equation}
which can be an approximate representation of $({\Omega_0^{\text{f}}},t=0)$, and
\begin{equation}
\begin{aligned}
    \mathcal{S}({\Omega_t^{\text{f}}}; [0,T]) = \{ (\boldsymbol{x}^{(i)}_{t^{(i)}},t^{(i)}) \quad|\quad &\boldsymbol{x}^{(i)}_{t^{(i)}}:= \boldsymbol{x}^{(i)}_{0} + \boldsymbol{d}^*(\boldsymbol{x}^{(i)}_{0},t^{(i)}; \boldsymbol{\theta}_{d}) \\
    &\text{ where }  \boldsymbol{x}^{(i)}_{0}\in\mathcal{S}({\Omega_0^{\text{f}}}) \text{ and }  t^{(i)}\in\mathcal{T}([0,T]) \}    
\end{aligned}
\label{eqnorm:ttt}
\end{equation}
which can be an approximate representation of $({\Omega_t^{\text{f}}},[0,T])$.

Naturally, we define the norm with $\mathcal{S}({\Omega_t^{\text{f}}}; [0,T])$ as
\begin{equation}
    ||f||^2_{\mathcal{S}({\Omega_t^{\text{f}}}; [0,T])} := \frac{1}{M} \sum\limits_{(\boldsymbol{x},t) \in \mathcal{S}({\Omega_t^{\text{f}}}; [0,T])} |f(\boldsymbol{x},t)|^2
\end{equation}
Similarly, we sample $N$ points for $\Gamma_0^{\text{w}}$ and $K$ points for $\Gamma_0^{\text{in}} \cup \Gamma_0^{\text{out}}$. Then $\mathcal{S}({\Gamma_t}; [0,T])$, $\mathcal{S}({\Gamma_0^{\text{w}}}; [0,T])$, $\mathcal{S}({\Omega_0^{\text{f}}}; [0,T])$, $\mathcal{S}({\partial\Gamma_0^{\text{w}}}; [0,T])$ can be obtained in the same way as \eqref{eqnorm:ttt}, while $\mathcal{S}({\partial\Gamma_0^{\text{w}}}; t=0)$ can refer to \eqref{eqnorm:t0}.
Then, we can compute for the losses in \eqref{eqs:fluid} and \eqref{eqs:solid} by redefining those norms as
\begin{equation}
    ||\cdot||^2_{L^2(\cdot; \cdot)} := ||\cdot||^2_{\mathcal{S}(\cdot; \cdot)} 
\end{equation}

\section{Parallel Training}

Neural network training can be parallelized in various ways to improve efficiency and reduce training time. We focus on data parallelism, a popular technique for parallelizing the training of neural networks, especially in distributed computing environments.

Taking the solid problem as an example, the network model $N_a$ is first replicated on $P$ processing units (CPU or GPU), where each replica is denoted as $N_a^{(i)}$ with $i=1,\dots,P$ indexing the processing unit for $P$ units. Each replica processes a different batch of data independently. Specifically, the sampled point $\mathcal{N}$ is divided into $P$ groups $\mathcal{N}^{(i)}$, where $i=1,\dots,P$. The forward pass, involving the computation of predictions, is performed independently on each replica as $N_a^{(i)}(\mathcal{N}^{(i)})$. The loss is computed for each batch as $\mathcal{L}_{\text{solid}}({\boldsymbol{\theta}_{d}}^{(i)};\mathcal{N}^{(i)})$, where ${\boldsymbol{\theta}_{d}}^{(i)}$ is the weight parameters for replica $N_a^{(i)}$. Backpropagation is performed independently on each model replica to calculate the gradients of the loss with respect to the model parameters to acquire the corresponding gradients ${\boldsymbol{\theta}_{d}}^{(i)}$. Gradients computed on each replica are then averaged and used to update the shared model parameters:
\begin{equation}
\nabla_{{\boldsymbol{\theta}_{d}}} 
 \mathcal{L}_{\text{solid}}(\boldsymbol{\theta}_{\text{d}};\mathcal{N}) = \frac{1}{P} \sum\limits_{i=1}^{P} \nabla_{{\boldsymbol{\theta}_{\text{d}}^{(i)}}} \mathcal{L}_{\text{solid}}(\boldsymbol{\theta}_{\text{d}}^{(i)};\mathcal{N}^{(i)})
\end{equation}

\end{document}